\documentclass[a4paper,notitlepage, eqno]{article}%
\usepackage{amssymb}
\usepackage{graphicx}
\usepackage{amsmath}
\usepackage{amsthm}
\usepackage{amsfonts}%
\providecommand{\U}[1]{\protect\rule{.1in}{.1in}}
\theoremstyle{plain}
\newtheorem{theorem}{Theorem}[section]

\newtheorem{cor}[theorem]{Corollary}

\newtheorem{prop}[theorem]{Proposition}

\newtheorem{thma}{Theorem}

\newtheorem{cora}[thma]{Corollary}

\newtheorem{exaa}[thma]{Example}

\theoremstyle{definition}

\newtheorem{rem}[theorem]{Remark}
\numberwithin{equation}{section}
\numberwithin{theorem}{section}

\allowdisplaybreaks 
\ifx\pdfoutput\relax\let\pdfoutput=\undefined\fi
\newcount\msipdfoutput
\ifx\pdfoutput\undefined\else
\ifcase\pdfoutput\else
\ifx\paperwidth\undefined\else
\ifdim\paperheight=0pt\relax\else\pdfpageheight\paperheight\fi
\ifdim\paperwidth=0pt\relax\else\pdfpagewidth\paperwidth\fi
\fi\fi\fi
\begin{document}

\title{Past and recent contributions to indefinite sublinear elliptic problems
 }
\author{Uriel Kaufmann\thanks{FaMAF-CIEM (CONICET), Universidad Nacional de
C\'{o}rdoba, Medina Allende s/n, Ciudad Universitaria, 5000 C\'{o}rdoba,
Argentina. \textit{E-mail address: }kaufmann@mate.uncor.edu} , Humberto Ramos
Quoirin\thanks{CIEM-FaMAF, Universidad Nacional de C\'{o}rdoba, (5000)
C\'{o}rdoba, Argentina. \textit{E-mail address: }humbertorq@gmail.com} ,
Kenichiro Umezu\thanks{Department of Mathematics, Faculty of Education,
Ibaraki University, Mito 310-8512, Japan. \textit{E-mail address:
}kenichiro.umezu.math@vc.ibaraki.ac.jp}
\and \noindent}
\maketitle

\begin{abstract}
We review the \textit{indefinite sublinear} elliptic equation 
%\marginpar{modified}
$-\Delta u=a(x)u^{q}$ in a smooth bounded domain
$\Omega\subset\mathbb{R}^{N}$, with Dirichlet or Neumann homogeneous boundary
conditions. Here $0<q<1$ and $a$ is continuous and changes sign, in which case
the strong maximum principle does \textit{not} apply. As a consequence, the
set of nonnegative solutions of these problems has a rich structure, featuring
in particular both \textit{dead core} and/or \textit{positive} solutions.
%\marginpar{modified}
Overall, we are interested in sufficient and necessary conditions on $a$ and
$q$ for the existence of \textit{positive} solutions. We describe the main
results from the past decades, and combine it with our recent contributions.
The proofs are briefly sketched.

\end{abstract}

%-------------------- abstract ---------------------

%------------------------ body ------------------

\section{Introduction}

\label{sec:Intro}

Let $N\geq1$, $\Omega\subset\mathbb{R}^{N}$ be a smooth bounded domain, and
$\Delta$ the usual Laplace operator. This article is devoted to the semilinear
equation
\begin{equation}
\label{e}-\Delta u=a(x)u^{q}\quad\mbox{in}\quad\Omega,
\end{equation}
under the condition
%\marginpar{put number}%
\[
a\mbox{ changes sign}\mbox{ and }0<q<1.\leqno{({\bf AQ})}
\]
This is a prototype of \textit{indefinite} (due to the change of sign of $a$)
and sublinear (with respect to $u$) elliptic pde, which is motivated by the
porous medium type equation \cite{GM77, PMeq}%
\[
w_{t}=\Delta(w^{m})+a(x)w,\quad m>1,
\]
after the change of variables $u=w^{m}$ and $q=1/m$. Indefinite elliptic
problems have attracted considerable attention since the 70's, mostly in the
linear ($q=1$) and superlinear ($q>1$) cases \cite{AT,ALG, BCN,BL,
fle,HK,LG,MM, Ou,Te}.
%\marginpar{ref 26 added}
We intend here to give an overview of the main results known in the
\textit{sublinear} case. For the sign-definite case $a\geq0$ we refer to
\cite{Am76, BO86, BK92, L73, L82}.

We shall consider \eqref{e} under Dirichlet and
%\marginpar{modified}
Neumann homogeneous boundary conditions, i.e. the problems
\[%
\begin{cases}
-\Delta u=a(x)u^{q} & \mbox{in $\Omega$},\\
u\geq0 & \mbox{in $\Omega$},\\
u=0 & \mbox{on $\partial \Omega$},
\end{cases}
\leqno{(P_{\mathcal{D}})}
\]
and
\[%
\begin{cases}
-\Delta u=a(x)u^{q} & \mbox{in $\Omega$},\\
u\geq0 & \mbox{in $\Omega$},\\
\partial_{\nu}u=0 & \mbox{on $\partial \Omega$},
\end{cases}
\leqno{(P_{\mathcal{N}})}
\]
where $\partial_{\nu}$ is the exterior normal derivative.

Throughout this article, we assume that $a\in C(\overline{\Omega})$. By a
\textit{solution} of $(P_{\mathcal{D}})$ we mean a \textit{strong}
%\marginpar{modified}
solution $u\in W_{\mathcal{D}}^{2,r}(\Omega)$ for some $r>N$, where
\[
W_{\mathcal{D}}^{2,r}(\Omega):=\{u\in W^{2,r}(\Omega):u=0\ \mbox{on}\ \partial
\Omega\}.
\]
Note that $u\in C^{1}(\overline{\Omega})$, and so the boundary condition is
satisfied in the usual sense. A similar definition holds for $(P_{\mathcal{N}%
})$. We say that a solution $u$ is \textit{nontrivial} if $u\not \equiv 0$,
and \textit{positive} if $u>0$ in $\Omega$. Among positive solutions of
$(P_{\mathcal{D}})$, we are interested in \textit{strongly positive} solutions
(denoted by $u\gg0$), namely, solutions in
\[
\mathcal{P}_{\mathcal{D}}^{\circ}:=\left\{  u\in C_{0}^{1}(\overline{\Omega
}):u>0\text{ in }\Omega\text{, \ and }\partial_{\nu}u<0\text{ on }%
\partial\Omega\right\}  .
\]
For $(P_{\mathcal{N}})$, a solution is \textit{strongly positive} if it
belongs to
\[
\mathcal{P}_{\mathcal{N}}^{\circ}:=\left\{  u\in C^{1}(\overline{\Omega
}):u>0\text{ on }\overline{\Omega}\right\}  .
\]
In case that \textit{every} nontrivial solution of $(P_{\mathcal{D}})$
(respect. $(P_{\mathcal{N}})$) is \textit{strongly positive} we say that this
problem has the \textit{positivity property}.

The condition (AQ) gives rise to the main feature of this class of problems,
namely, the fact that the strong maximum principle (shortly \textbf{SMP}) does
\textit{not} apply. Let us recall the following version of this result
%\marginpar{modified}
(for a proof, see e.g. \cite[Theorem 7.10]{LGbook13}):\newline

\noindent\textbf{Strong maximum principle:} Let $u\in W^{2,r}(\Omega)$
%\marginpar{changed to Thm. A}
for some $r>N$ be such that $u\geq0$ and $(-\Delta+M)u\geq0$ in $\Omega$, for
some constant $M\geq0$. Then either $u\equiv0$ or $u>0$ in $\Omega$ and
$\partial_{\nu}u(x)<0$ for any $x\in\partial\Omega$ such that $u(x)=0$%
.\newline

Given $u$ satisfying \eqref{e},
%\marginpar{modified}
we see that under (AQ) we can't find in general some $M>0$ such that
$(-\Delta+M)u=a(x)u^{q}+Mu\geq0$ in $\Omega,$ which prevents us to apply the
\textbf{SMP}, unlike when $a\geq0$ (the \textit{definite} case) or $q\geq1$
(the \textit{linear} and \textit{superlinear} cases). This fact is reinforced
by a simple example of a nontrivial solution $u$ (of both $(P_{\mathcal{D}})$
and $(P_{\mathcal{N}})$) violating the conclusion of the \textbf{SMP} (see
Example \ref{exa:cos} below), which shows that the positivity property may
fail. Moreover, such example also provides us with nontrivial \textit{dead
core} solutions of $(P_{\mathcal{D}})$ and $(P_{\mathcal{N}})$, i.e. solutions
vanishing in some open subset of $\Omega$.
%\newline
%(e.g. \cite[Remark 3.7]{KRQU2019}).

To the best of our knowledge, the study of $\left(  P_{\mathcal{D}}\right)  $
and $\left(  P_{\mathcal{N}}\right)  $ in the indefinite and sublinear case
was launched in the late 80's by Bandle, Pozio and Tesei \cite{BPT1,BPT2,PT}.
These works were then followed by the contributions of Hern\'{a}ndez, Mancebo
and Vega \cite{HMV}, Delgado and Suarez \cite{DS}, and Godoy and Kaufmann
\cite{GK1,GK2}. We shall review the main results of these papers in the next
section and complement it with our main recent results from \cite{KRQU16,
KRQUnodea, KRQU3} in the subsequent sections. Since the proofs can be found in
the aforementioned articles,
%\marginpar{modified, to avoid repetition}
in most cases we shall only sketch them here.

\section{First results}

Let us recall the first existence and uniqueness results on the problems
above. For the Neumann problem, the following condition on $a$ plays an
important role:
\[
\int_{\Omega}a<0.\leqno{(\mbox{\bf A.0})}
\]
Indeed,
%\marginpar{\textit{modified and added}}
we shall see that (A.0) is \textit{necessary} for the existence of a positive 
%\marginpar{changed 'positive' by 'nontrivial', since at this point we're recalling the results from Bandle et al}
solution of $(P_{\mathcal{N}})$, and
\textit{sufficient} for the existence of a nontrivial solution, for any
$q\in(0,1)$. As for the uniqueness results, some merely technical conditions
(see also the beginning of Section \ref{sec:oq}) on the set
\[
\Omega_{+}:=\{x\in\Omega:a(x)>0\}
\]
shall be used, namely:
\[
\Omega_{+}\text{ has \textit{finitely }many connected components,}%
\leqno{(\mbox{\bf A.1})}
\]%
\[
\partial\Omega_{+}\text{ satisfies an inner sphere condition with respect to
}\Omega_{+}.\leqno{(\mbox{\bf A.2})}
\]

The following results were proved by Bandle, Pozio and Tesei \cite{BPT1,BPT2},
and Delgado and Su\'{a}rez \cite{DS}. Although \cite{BPT1,BPT2} require that
$a\in C^{\theta}(\overline{\Omega})$ for some $0<\theta<1$, one can easily see
from the proofs that these results still hold for strong solutions assuming
that $a\in C(\overline{\Omega})$.

\begin{thma}
\strut\label{t10}

\begin{enumerate}
\item The Dirichlet case:

\begin{enumerate}
\item $(P_{\mathcal{D}})$ has at most one positive solution \cite[Theorem
2.1]{DS}. Moreover, if (A.1) and (A.2) hold then $(P_{\mathcal{D}})$ has at
most one solution positive in $\Omega_{+}$ \cite[Theorem 2.1]{BPT1}.

\item $(P_{\mathcal{D}})$ has at least one nontrivial solution \cite[Theorem
2.2]{BPT1}.
\end{enumerate}

\item The Neumann case:

\begin{enumerate}
\item $(P_{\mathcal{N}})$ has at most one solution in $\mathcal{P}
_{\mathcal{N}}^{\circ}$ \cite[Lemma 3.1]{BPT2}. Moreover, if (A.1) and (A.2)
hold then $(P_{\mathcal{N}})$ has at most one solution positive in $\Omega
_{+}$ \cite[Theorem 3.1]{BPT2}.

\item If (A.0) holds then $(P_{\mathcal{N}})$ has at least one nontrivial
solution. Conversely, if $(P_{\mathcal{N}})$ has a positive solution then
(A.0) holds \cite[Theorem 2.1]{BPT2}.
\end{enumerate}
\end{enumerate}
\end{thma}

%\begin{proof}
\textit{Sketch of the proof}. The uniqueness assertions rely on the following
change of variables: if $u>0$ and $-\Delta u=a(x)u^{q}$ in $\Omega$ then
$v:=(1-q)^{-1}u^{1-q}$ solves $-\Delta v=qu^{q-1}|\nabla v|^{2}+a(x)$ in
$\Omega$. Let $u_{1},u_{2}$ be positive solutions of $(P_{D})$ and assume that
$\tilde{\Omega}:=\{x\in\Omega:u_{1}(x)>u_{2}(x)\}$ is nonempty. We set
$v_{i}:=(1-q)^{-1}u_{i}^{1-q}$ for $i=1,2$, so that $\Phi:=v_{1}-v_{2}>0$ in
$\tilde{\Omega}$. In addition,
\[
-\Delta\Phi=q\left(  u_{1}^{q-1}|\nabla v_{1}|^{2}-u_{2}^{q-1}|\nabla
v_{2}|^{2}\right)  <qu_{1}^{q-1}\left(  |\nabla v_{1}|^{2}-|\nabla v_{2}%
|^{2}\right)  ,
\]
i.e.
\begin{equation}
-\Delta\Phi-qu_{1}^{q-1}\nabla(v_{1}+v_{2})\nabla\Phi<0\quad\mbox{in}\quad
\tilde{\Omega}. \label{ef}%
\end{equation}
Since $\Phi=0$ on $\partial\tilde{\Omega}$, we obtain a contradiction with the
maximum principle. This shows that $(P_{\mathcal{D}})$ has at most one
positive solution. Now, if $u_{1},u_{2}\in\mathcal{P}_{\mathcal{N}}^{\circ}$
solve $(P_{\mathcal{N}})$ then $\Phi$ satisfies \eqref{ef} and for any
$x\in\partial\tilde{\Omega}$ we have either $\Phi(x)=0$ or $\partial_{\nu}%
\Phi(x)=0$. By the maximum principle, we infer that $\Phi$ is constant in
$\tilde{\Omega}$, which contradicts \eqref{ef}. The proof of the uniqueness of
a solution of $(P_{\mathcal{D}})$ positive in $\Omega_{+}$ (respect. a
solution of $(P_{\mathcal{N}})$ in $\mathcal{P}_{\mathcal{N}}^{\circ}$) uses
the same change of variables, but is more involved. We refer to
\cite{BPT1,BPT2} for the details.

The existence results can be proved either by a variational argument or by the
sub-supersolutions method. In the first case, it suffices to show that the
functional
\[
I_{q}(u):=\int_{\Omega}\left(  \frac{1}{2}|\nabla u|^{2}-\frac{1}%
{q+1}a(x)|u|^{q+1}\right)
\]
has a negative global minimum in $H_{0}^{1}(\Omega)$ or $H^{1}(\Omega)$. In
the latter case the condition (A.0) is crucial. The second approach consists
in taking a ball $B\subset\Omega_{+}$ and a sufficiently small first positive
eigenfunction of $-\Delta$ on $H_{0}^{1}(B)$ extended by zero to $\Omega$, to
find a (nontrivial) subsolution of both $(P_{\mathcal{D}})$ and
$(P_{\mathcal{N}})$. An arbitrary large supersolution of $(P_{\mathcal{D}})$
is given by $kz$, where $z$ is the unique solution of $-\Delta z=a^{+}$ in
$\Omega$, $z=0$ on $\partial\Omega$, and $k>0$ is large enough (as usual, we
write $a=a^{+}-a^{-}$, with $a^{\pm}:=\max(\pm a,0)$). The construction of 
a 
%\marginpar{replaced 'big' by 'suitable'}
suitable supersolution of
$(P_{\mathcal{N}})$ under (A.0) is more delicate, and we refer to \cite{BPT2}
for the details.

Finally, if $(P_{\mathcal{N}})$ has a positive solution $u$ then, multiplying
the equation by $(u+\varepsilon)^{-q}$ (with $0<\varepsilon<1$) and
integrating by parts, we find that
\[
\int_{\Omega}a\left(  \frac{u}{u+\varepsilon}\right)  ^{q}=-q\int_{\Omega
}(u+\varepsilon)^{-(q+1)}|\nabla u|^{2}<-q\int_{\Omega}(u+1)^{-(q+1)}|\nabla
u|^{2}<0.
\]
Letting $\varepsilon\rightarrow0$ we can check that $\int_{\Omega}a<0$.
\qed\newline

Although not stated explicitly in \cite{BPT1,BPT2}, the next corollary follows
almost directly from the existence and uniqueness results in these papers.

\begin{cora}
\label{A=I} Let $\Omega_{+}$ be connected and satisfy (A.2). Then
$(P_{\mathcal{D}})$ has a unique nontrivial solution. The same conclusion
holds for $(P_{\mathcal{N}})$ assuming in addition (A.0).
\end{cora}

\textit{Sketch of the proof}. It is based on the fact that a nontrivial
solution $u$ of $(P_{\mathcal{D}})$ or $(P_{\mathcal{N}})$ satisfies
$u\not \equiv 0$ in $\Omega_{+}$, which follows from the inequality
$0<\int_{\Omega}|\nabla u|^{2}\leq\int_{\Omega}a^{+}(x)u^{q+1}$.
%\marginpar{removed}
Since $\Omega_{+}$ is connected, by the maximum principle we find that $u>0$
in $\Omega_{+}$. And there is only one solution having this property, by
Theorem \ref{t10}. \qed\newline

\begin{rem}
\strut

\begin{enumerate}
\item Let us remark that the nontrivial solutions provided by Theorem
\ref{t10} (i-b) and (ii-b) are \textit{not} necessarily unique, see e.g.
\cite{BPT1,BPT2}.

\item Regarding Theorem \ref{t10} (ii-b), it is worth pointing out that (A.0)
is \textit{not} necessary for the existence of a nontrivial solution of
$(P_{\mathcal{N}})$ for some $q\in(0,1)$, cf. \cite[Section 4]{BPT2} and
\cite[Remark 4.3]{KRQUnodea}.
\end{enumerate}
\end{rem}

\begin{figure}[tbh]
\centerline{
		\includegraphics[scale=0.125]{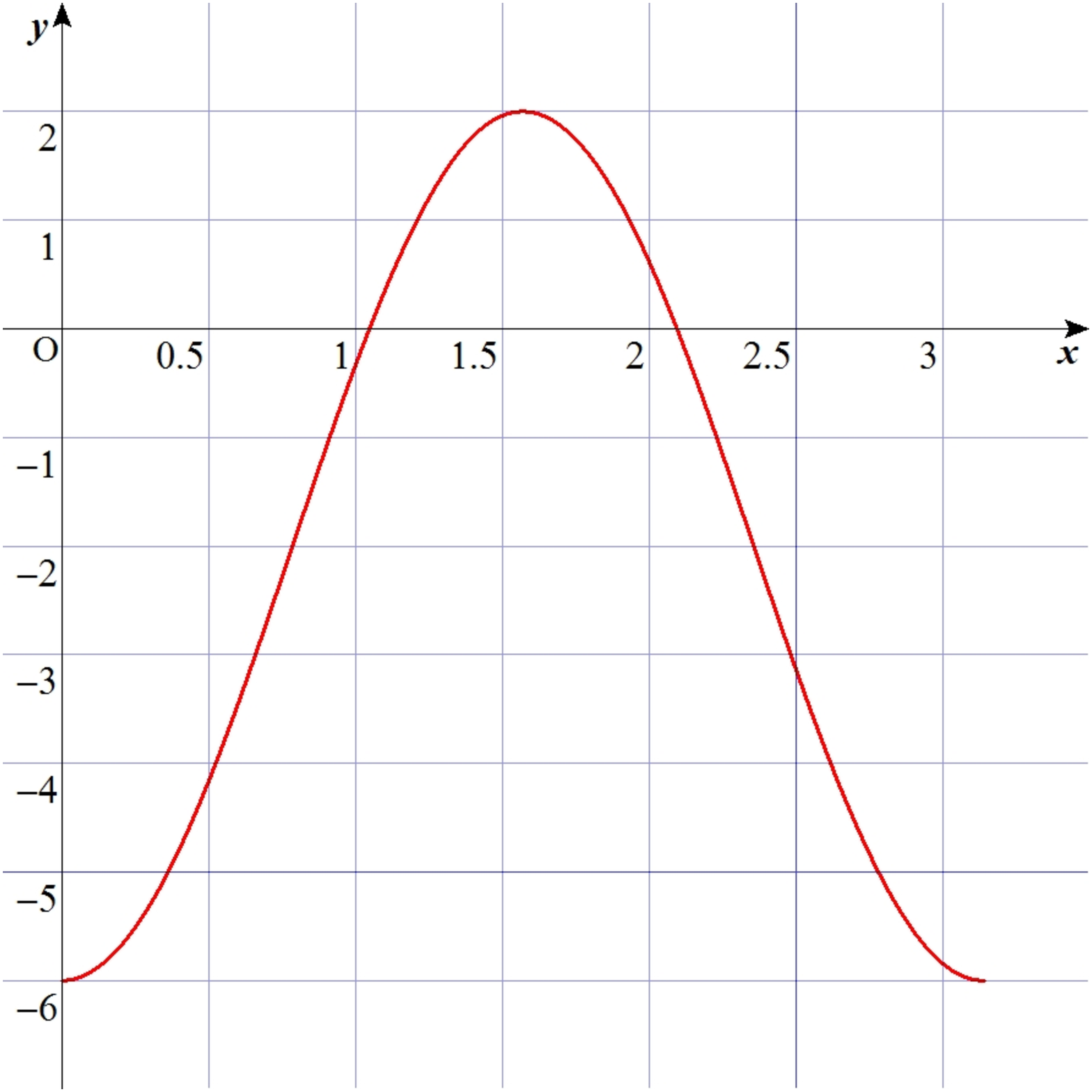}
		\hskip0.01cm
		\includegraphics[scale=0.125]{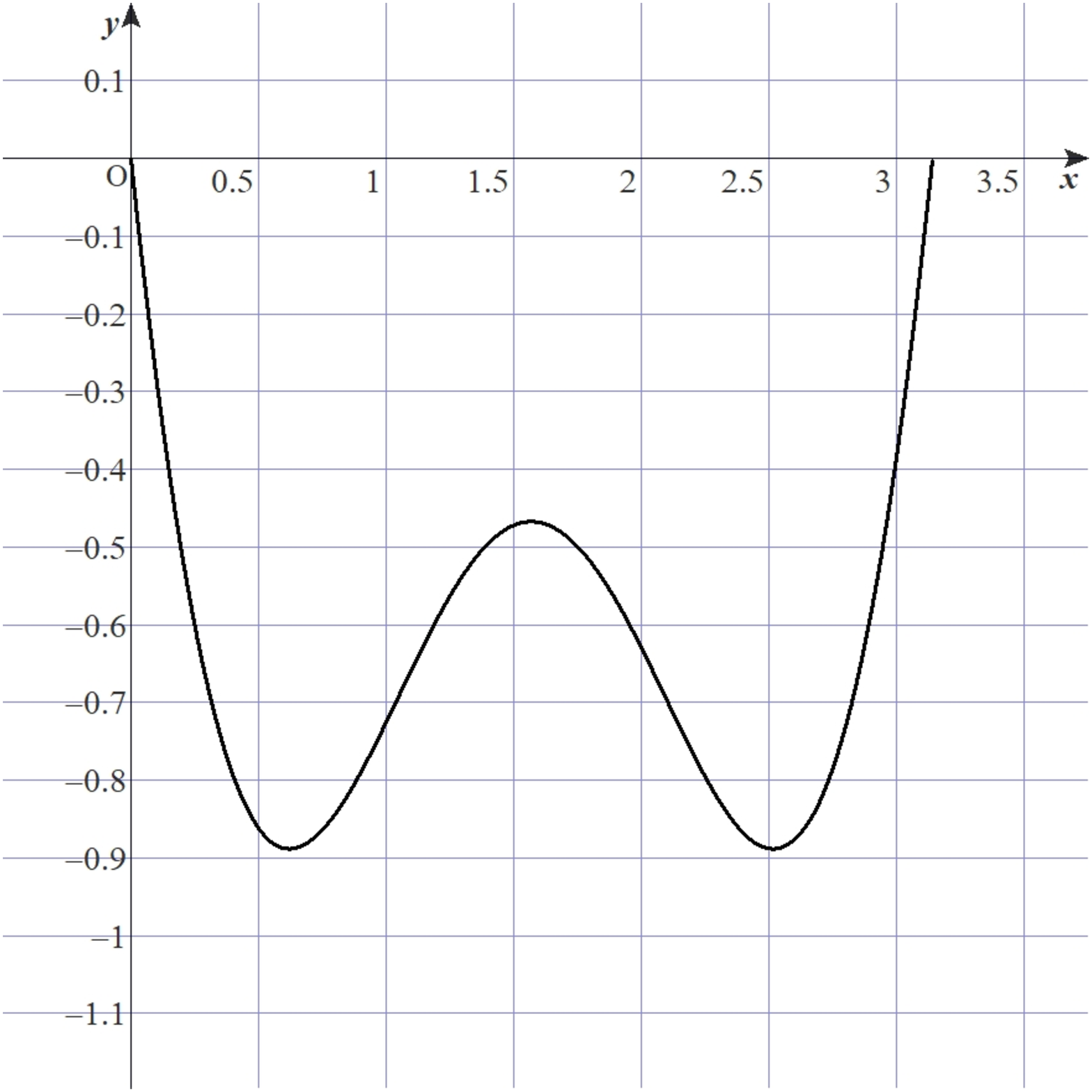}
		\hskip0.01cm
		\includegraphics[scale=0.125]{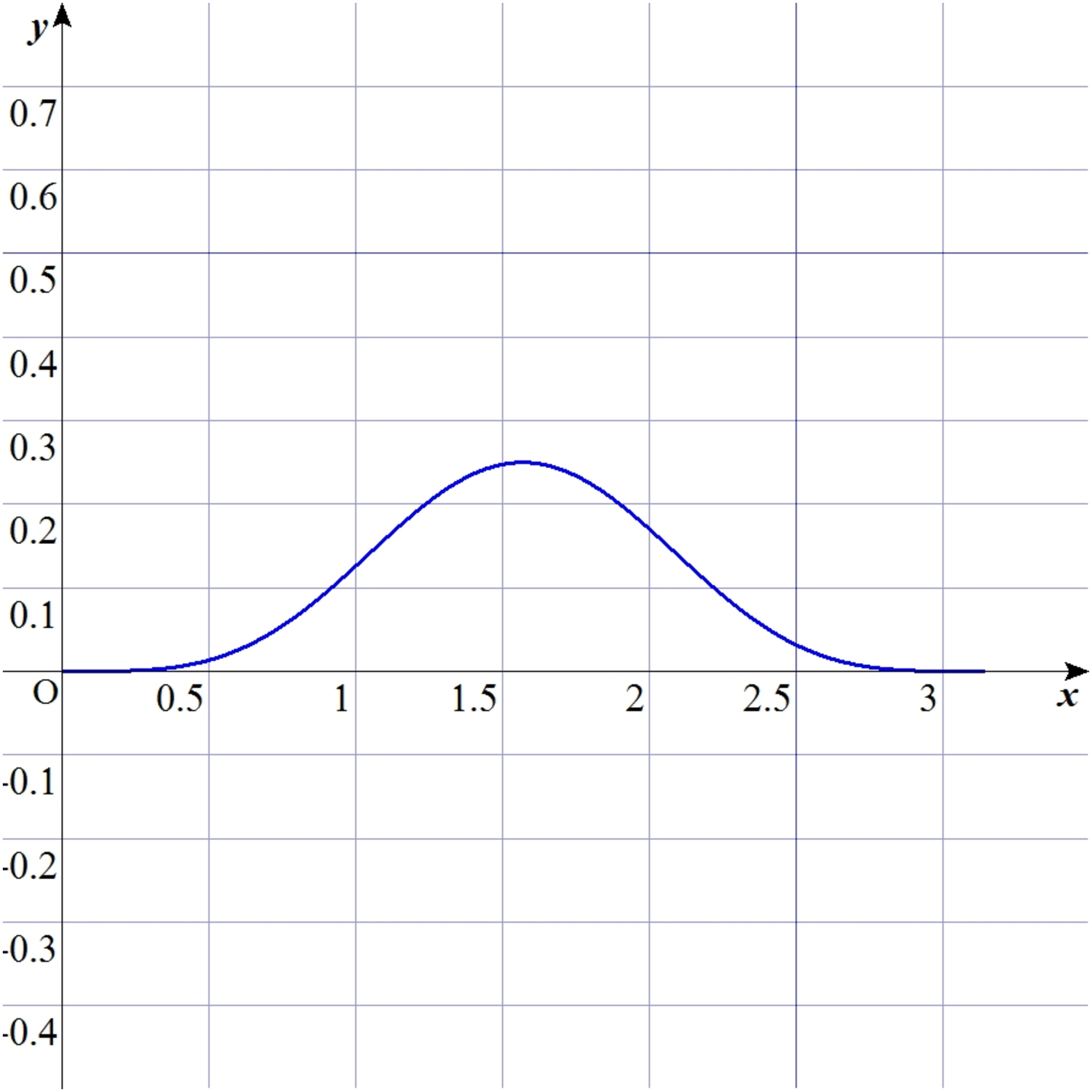}
	} \centerline{(i) \hskip3.55cm 	(ii)\hskip3.55cm (iii) }\caption{(i) The
indefinite weight $a_{\frac{1}{2}}$; (ii) $\mathcal{S}(a_{\frac{1}{2}})$;
(iii) The positive solution $u\not \gg 0$ for $a_{\frac{1}{2}}$.}%
\label{figcexa}%
\end{figure}

Let us now
%\marginpar{added}
give an example of a nontrivial solution $u\not \gg 0$ of $(P_{\mathcal{D}})$
and $(P_{\mathcal{N}})$. It is essentially due to \cite{GK1}, where the case
$q=\frac{1}{2}$ was considered (see Figure \ref{figcexa}).
%\marginpar{\textbf{added}}

\begin{exaa}
\label{exa:cos}

Let $\Omega:=(0,\pi)$ and $q\in(0,1)$. We choose
\[
r=r_{q}:=\frac{2}{1-q}\in\left(  2,\infty\right)  ,\quad a(x)=a_{q}%
(x):=r^{1-\frac{2}{r}}\left(  1-r\cos^{2}x\right)  \quad\text{for }%
x\in\overline{\Omega}.
\]
Then $u(x):=\frac{\sin^{r}x}{r}\in C^{2}(\overline{\Omega})$ satisfies
\[%
\begin{cases}
-u^{\prime\prime}=a(x)u^{q} & \mbox{in $\Omega$},\\
u>0 & \mbox{in $\Omega$},\\
u=u^{\prime}=u^{\prime\prime}=0 & \mbox{on $\partial \Omega$}.
\end{cases}
\]

\end{exaa}

The above example also provides \textit{dead core} solutions of
$(P_{\mathcal{D}})$ and $(P_{\mathcal{N}})$ both. Indeed, it suffices to
consider any bounded open interval $\Omega^{\prime}$ with $\Omega^{\prime
}\supset\overline{\Omega}$, and extend $u$ by zero and $a$ in any way to
$\Omega^{\prime}$. Then $u$ is a nontrivial dead core solution of both
$(P_{\mathcal{D}})$ and $(P_{\mathcal{N}})$, considered now in $\Omega
^{\prime}$.

Since the \textbf{SMP} does not apply and
%\marginpar{modified}
\textit{dead core} solutions may exist, obtaining a positive solution for
these problems is a delicate issue which has been given little consideration.
%\marginpar{modified}
Let $\varphi\in W_{\mathcal{D}}^{2,r}\left(  \Omega\right)  $ be the unique
solution of the Poisson equation
\[%
\begin{cases}
-\Delta\varphi=a(x) & \mbox{ in }\Omega,\\
\varphi=0 & \mbox{ on }\partial\Omega,
\end{cases}
\]
and $\mathcal{S}:L^{r}(\Omega)\rightarrow W_{\mathcal{D}}^{2,r}(\Omega)$ be
the corresponding \textit{solution operator}, i.e. $\mathcal{S}(a)=\varphi$.
In \cite{HMV} Hern\'{a}ndez, Mancebo and Vega showed that the condition
\[
\mathcal{S}(a)\gg0\leqno{(\mbox{\bf A.3})}
\]
%\marginpar{sentence removed}
implies the existence of a positive solution of $(P_{\mathcal{D}})$ for all
$q\in\left(  0,1\right)  $. Later on Godoy and Kaufmann \cite{GK1,GK2}
provided other sufficient conditions, namely, that $a^{-}$ is small enough, or
$q$ is close enough to $1$ (for some particular choices of $N$, $\Omega$, and
$a$). We shall state a simplified version of these results in the sequel, and
refer to \cite[Theorem 4.4]{HMV}, \cite[Theorems 3.1 and 3.2]{GK2}, and
\cite[Theorems 2.1 (i) and 3.2]{GK1} for the precise statements.
%\marginpar{moved sentence earlier}

\begin{thma}
\strut\label{t11}

\begin{enumerate}
\item If $a$ satisfies (A.3) then $(P_{\mathcal{D}})$ has a positive solution
for every $q \in(0,1)$.

\item Let $q$ and $a^{+}$ be fixed. Then there exists a constant $C>0$ such
that $(P_{\mathcal{D}})$ has a positive solution if $\|a^{-}\|_{C(\overline
{\Omega})}<C$.

\item If either $N=1$ or $\Omega$ is a ball, $a$ is radial, and $0\not \equiv
a\geq0$ in some smaller ball, then there exists $\overline{q}=\overline{q}(a)$
such that $(P_{\mathcal{D}})$ has a positive solution for $\overline{q}%
<q<1$.\newline
\end{enumerate}
\end{thma}

\begin{rem}
Let us mention that Theorem \ref{t11} (i) is still true for a linear second
order elliptic operator with nonnegative zero order coefficient. On the other
side, it may happen that $\mathcal{S}(a)<0$ \textit{everywhere} in $\Omega$
and yet $(P_{\mathcal{D}})$ admits a positive solution for some $q\in\left(
0,1\right)  $. Indeed, if we take $q=\frac{1}{2}$ in Example \ref{exa:cos}
then $\mathcal{S}(a_{q})=x^{2}-\pi x+1-\cos2x<0$ in $\left(  0,\pi\right)  $,
see Figure \ref{figcexa} (ii). Note also that (A.3) is \textit{not} compatible
with the existence of a positive solution
%\marginpar{modified}
for $(P_{\mathcal{N}})$, since it implies $\int_{\Omega}a>0$,
%\marginpar{\textit{modified}}
contradicting (A.0), which is necessary by Theorem \ref{t10} (ii-b).
\end{rem}

%\begin{figure}[tbh]
%\begin{center}
%\includegraphics[scale=0.17]{fig_sa2.eps}
%\end{center}
%\caption{$\mathcal{S}(a_{\frac{1}{2}})$ in Example \ref{exa:cos}. }%
%\label{figIntro}%
%\end{figure}

\textit{Sketch of the proof}. All assertions follow by the well known
sub-supersolutions method. Let us note that (\textit{unlike} for
$(P_{\mathcal{N}})$) it is easy to provide arbitrary big supersolutions for
$(P_{\mathcal{D}})$. Indeed, a few computations show that $k\mathcal{S}\left(
a^{+}\right)  $ is a supersolution of $(P_{\mathcal{D}})$ for all $k>0$ large
enough. So the only task is to provide a \textit{positive }subsolution. In
(i), after some computations one can check that $\left[  \left(  1-q\right)
\mathcal{S}\left(  a\right)  \right]  ^{1/\left(  1-q\right)  }$ is the
desired subsolution.

In both (ii) and (iii), the subsolution is constructed by splitting the domain
in two parts (a ball $B$ in which $0\not \equiv a\geq0$, and $\Omega\setminus
B$), constructing \textquotedblleft subsolutions\textquotedblright\ in each of
them, and checking that they can be glued appropriately to get a subsolution
in the entire domain (see \cite{bl}). This fact depends on obtaining estimates
for the normal derivatives of these subsolutions on the boundaries of the
subdomains. In (iii) these bounds can be computed rather explicitly using the
symmetry of $a$ and the fact that $\Omega$ is a ball, while in (ii) the key
tool is an estimate due to Morel and Oswald \cite[Lemma 2.1]{cabre}. The proof
of both (ii) and (iii) involve several computations, and we refer to
\cite{GK1,GK2} for the details. \qed\newline

%-----
%-------

Godoy and Kaufmann \cite{GK2} also proved that when $a$ is too negative in a
ball there are no positive solutions of $(P_{\mathcal{D}})$ (see also Remark
\ref{r} (i)
%\marginpar{modified}
below). This result can also be seen as a first step towards the construction
of dead core solutions.

\begin{thma}
\label{t12}

Let $q$ and $a^{+}$ be fixed. Given a ball $B=B_{R}(x_{0})\subset
\Omega\setminus\Omega_{+}$ there exists a constant $C=C(\Omega,N,q,R,a^{+})>0$
such that any solution of $(P_{\mathcal{D}})$ vanishes at $x_{0}$ if
$\min_{\overline{B}}a^{-}>C$.
\end{thma}

\textit{Sketch of the proof}. We use a comparison argument: let $u$ be a
nontrivial solution of $(P_{\mathcal{D}})$, and $\underline{a}:=\min
_{\overline{B}}a^{-}$. Set
\[
C_{N,q}:=\frac{\left(  1-q\right)  ^{2}}{2\left(  N\left(  1-q\right)
+2q\right)  }\quad\text{and\quad}w(x):=\left(  C_{N,q}\underline{a}\left\vert
x-x_{0}\right\vert ^{2}\right)  ^{\frac{1}{1-q}}.
\]
One can check that $\Delta w\leq a^{-}w^{q}$ in $B$. On the other hand, note
that $\Delta u=a^{-}u^{q}$ in $B$ and $\left\Vert u\right\Vert _{\infty}%
\leq\left(  \left\Vert \mathcal{S}\right\Vert \left\Vert a^{+}\right\Vert
_{\infty}\right)  ^{{\frac{1}{1-q}}}$, and so $u\leq w$ on $\partial B$ if%
\begin{equation}
\underline{a}\geq\frac{\left\Vert \mathcal{S}\right\Vert \left\Vert
a^{+}\right\Vert _{\infty}}{R^{2}C_{N,q}}. \label{ia}%
\end{equation}
It follows then from the comparison principle that $u\leq w$ in $B$. In
particular, $u\left(  x_{0}\right)  =0$. \qed

\begin{rem}
\strut\label{r}

\begin{enumerate}
\item The latter proof can be adapted for the Neumann problem, taking into
account the following \textit{a priori} bound: Under (A.0), there exists $C>0$
(independent of $a^{-}$) such that
%\marginpar{modified, since in Nodea we don't mention that the constant does not depend on $a^-$}
$\Vert u\Vert_{C(\overline{\Omega})}\leq C$ for every subsolution of
$(P_{\mathcal{N}})$.
%\marginpar{removed ref. to Nodea}

\item Note that $C_{N,q}\rightarrow0$ as $q\rightarrow1^{-}$, i.e. the closer
is $q$ to $1$, the larger is the right-hand side in \eqref{ia}, and the more
negative $a$ needs to be in $B_{R}(x_{0})$ to satisfy \eqref{ia}. This fact is
consistent with Theorem \ref{t11} (ii) and (iii).
\end{enumerate}
\end{rem}

\section{Recent results}

Let us now briefly describe our main contributions to the study of
$(P_{\mathcal{D}})$ and $(P_{\mathcal{N}})$, which can be found in
\cite{KRQU16, KRQUnodea, KRQU3}:

\begin{enumerate}
\item[(I)] We determine the values of $q\in(0,1)$ for which $(P_{\mathcal{D}%
})$ and $(P_{\mathcal{N}})$ have the positivity property. In other words, we
provide a characterization of the following \textit{positivity sets}:
\begin{align*}
&  \mathcal{A}_{\mathcal{D}}=\mathcal{A}_{\mathcal{D}}(a):=\{q\in
(0,1):\mbox{$u\gg 0$ for any nontrivial solution $u$ of $(P_{\mathcal{D}})$}\},\\
&  \mathcal{A}_{\mathcal{N}}=\mathcal{A}_{\mathcal{N}}(a):=\{q\in
(0,1):\mbox{$u\gg 0$ for any nontrivial solution $u$ of $(P_{\mathcal{N}})$}\}.
\end{align*}
Thanks to a continuity argument inspired by Jeanjean \cite{Je}, and based on
the fact that the \textbf{SMP} applies when $q=1$, we shall see
%\marginpar{added}
in Theorem \ref{t2A} that under (A.1) we have $\mathcal{A}_{\mathcal{D}%
}=(q_{\mathcal{D}},1)$ and, assuming additionally (A.0), $\mathcal{A}%
_{\mathcal{N}}=(q_{\mathcal{N}},1)$, for some $q_{\mathcal{D}},q_{\mathcal{N}%
}\in\left[  0,1\right)  $ (see also Corollary \ref{old6.2}
%\marginpar{added}
and Theorem \ref{t3a}).

Note that in view of the existence and uniqueness results in Theorem
\ref{t10}, the sets $\mathcal{A}_{\mathcal{D}}$ and $\mathcal{A}_{\mathcal{N}%
}$ can also be expressed as follows:
\begin{align}
&  {\mathcal{A}}_{\mathcal{D}}=\{q\in(0,1):\left(  P_{\mathcal{D}}\right)
\text{ has a unique nontrivial solution }u\text{, and }u\gg0\},\label{3D}\\
&  {\mathcal{A}}_{\mathcal{N}}=\{q\in(0,1):\left(  P_{\mathcal{N}}\right)
\text{ has a unique nontrivial solution }u\text{, and }u\gg0\}.\nonumber
\end{align}

We also obtain some positivity properties for the \textit{ground state}
solution of $(P_{\mathcal{D}})$.
%\begin{enumerate}

\item[(II)] By the previous discussion we deduce that $(P_{\mathcal{D}})$
(respect.\ $(P_{\mathcal{N}})$, under (A.0)) has a solution $u\gg0$ for
$q\in\mathcal{A}_{\mathcal{D}}$ (respect.\ $q\in\mathcal{A}_{\mathcal{N}}$).
Thus, setting
\begin{align*}
&  \mathcal{I}_{\mathcal{D}}=\mathcal{I}_{\mathcal{D}}(a):=\left\{
q\in(0,1):\mbox{$(P_{\mathcal{D}})$ has a solution $u\gg 0$}\right\}  ,\\
&  \mathcal{I}_{\mathcal{N}}=\mathcal{I}_{\mathcal{N}}(a):=\left\{
q\in(0,1):\mbox{$(P_{\mathcal{N}})$ has a solution $u\gg 0$}\right\}  ,
\end{align*}
we observe that $\mathcal{A}_{\mathcal{D}}\subseteq\mathcal{I}_{\mathcal{D}}$
and $\mathcal{A}_{\mathcal{N}}\subseteq\mathcal{I}_{\mathcal{N}}$. We will
further investigate $\mathcal{I}_{\mathcal{D}}$ (respect. $\mathcal{I}%
_{\mathcal{N}}$) and analyze how close $\mathcal{A}_{\mathcal{D}}$ and
$\mathcal{I}_{\mathcal{D}}$ (respect.\ $\mathcal{A}_{\mathcal{N}}$ and
$\mathcal{I}_{\mathcal{N}}$) can be to each other,
%\marginpar{added}
see Theorems \ref{mt2c}, \ref{th}, and also Proposition \ref{ti} and Remark
\ref{remgs} (i).

%\marginpar{added}.
%\end{enumerate}
Note that Corollary \ref{A=I} tells us that
%\marginpar{\it\footnotesize is this corollary true for Dirichlet ? are the computations valid similarly as in Neumann ?}
if $\Omega_{+}$ is connected and satisfies (A.2), then $\mathcal{A}%
_{\mathcal{D}}=\mathcal{I}_{\mathcal{D}}$, and if additionally (A.0) holds,
then %\marginpar{modified}
$\mathcal{A}_{\mathcal{N}}=\mathcal{I}_{\mathcal{N}%
}$. Assuming moreover (A.3), we find by Theorem \ref{cormt2} (iv-c) that
$\mathcal{A}_{\mathcal{D}}=\left(  0,1\right)  $.
%\marginpar{modified sentence}
%\marginpar{added text}

%\begin{enumerate}

\item[(III)] We consider $(P_{\mathcal{D}})$ and $(P_{\mathcal{N}})$ via a
bifurcation approach, looking at
%\marginpar{modified}
$q$ as a bifurcation parameter and taking advantage of the fact that
$(P_{\mathcal{D}})$ has a trivial line of strongly positive solutions when
$q=1$, see Theorems \ref{cormt2}
%\marginpar{modified}
and \ref{c1} for $(P_{\mathcal{D}})$ and $(P_{\mathcal{N}})$, respectively.
%Proposition \ref{prop:gam2}.
We also analyze the structure of the nontrivial solutions set (with respect to
$q$) of $(P_{\mathcal{D}})$ and $(P_{\mathcal{N}})$ via variational methods
and the construction of sub and supersolutions, see Theorem
%\marginpar{modified}
\ref{cormt2} and Remark \ref{rem:nss} for $(P_{\mathcal{D}})$; Remarks
\ref{rem:geN} and \ref{rem:nss} for $(P_{\mathcal{N}})$. In particular, we
describe the asymptotic behaviors of nontrivial solutions as $q\rightarrow
0^{+}$ and $q\rightarrow1^{-}$.
%
%\textit{Comment: I think here we should try to be more specific and\marginpar{see comment} quote our main results in these issues, because otherwise I think the main results are difficult to find for the reader.}

\item[(IV)]
%\marginpar{item added}
Finally, in Section 6 we present, without proofs, two further kind of results.
On the one hand, we provide \textit{explicit }sufficient conditions for the
existence of positive solutions for $(P_{\mathcal{D}})$ and $(P_{\mathcal{N}%
})$, see Theorems \ref{cc} and \ref{rad2}. And 
%\marginpar{\textit{modified,because ?? appeared}}
on the other hand, in Theorem \ref{dc} we state
sufficient conditions for the existence of dead core solutions for
$(P_{\mathcal{N}})$.\smallskip
\end{enumerate}

The above issues will be developed in the forthcoming sections. In the last
section we include some final remarks and list some open questions.

\section{The positivity property}

\label{sec:P}

The next theorem extends Theorem \ref{t11} (iii)
%\marginpar{\it\footnotesize added}
under (A.1), showing that $(P_{\mathcal{D}})$, as well as $(P_{\mathcal{N}})$
under $(A.0)$, has a positive solution (and no other nontrivial solution) if
$q$ is close enough to $1$. In other words, we show that under (A.1) the
positivity property holds for such values of $q$ \cite[Theorems 1.3 and
1.7]{KRQU16}:
%\marginpar{modified}

\begin{theorem}
\label{t2A}Assume (A.1). Then:

\begin{enumerate}
\setlength{\itemsep}{0.2cm}

\item[(i)] $\mathcal{A}_{\mathcal{D}}=(q_{\mathcal{D}},1)$ for some
$q_{\mathcal{D}} \in[0,1)$.
%\marginpar{removed the second statements, since these ones are stated in Rem 1.3}

\item[(ii)] If (A.0) holds then $\mathcal{A}_{\mathcal{N}}=(q_{\mathcal{N}
},1)$ for some $q_{\mathcal{N}} \in[0,1)$.
\end{enumerate}
\end{theorem}

\textit{Sketch of the proof}. First we show that $\mathcal{A}_{\mathcal{D}}$
is nonempty. We proceed by contradiction, assuming that $q_{n}\rightarrow
1^{-}$ and $u_{n}$ are nontrivial solutions of $(P_{\mathcal{D}})$ with
$q=q_{n}$ and $u_{n}\not \gg 0$. We know that $u_{n}\not \equiv 0$ in
$\Omega_{+}$, and thanks to (A.1) we can assume that, for every $n\in
\mathbb{N}$, $u_{n}>0$ in some fixed connected component of $\Omega_{+}$. If
$\{u_{n}\}$ is bounded in $H_{0}^{1}(\Omega)$ then, by standard compactness
arguments, up to a subsequence, we have $u_{n}\rightarrow u_{0}$ in $H_{0}%
^{1}(\Omega)$ and $u_{0}$ solves $-\Delta u_{0}=a(x)u_{0}$. Moreover, we can
show that $\{u_{n}\}$ is away from zero, so that $u_{0}\not \equiv 0$. By the
\textbf{SMP} we get that $u_{0}\gg0$. Finally, by standard elliptic
regularity, we find that $u_{n}\rightarrow u_{0}$ in $C^{1}(\overline{\Omega
})$, up to a subsequence. Thus $u_{n}\gg0$ for $n$ large enough, and we have a
contradiction. If $\{u_{n}\}$ is unbounded in $H_{0}^{1}(\Omega)$ then,
normalizing it, we obtain a sequence $v_{n}$ converging to some $v_{0}%
\not \equiv 0$ that solves an eigenvalue problem. Once again, the \textbf{SMP}
implies that $v_{0}\gg0$, a contradiction. A similar argument shows that
$\mathcal{A}_{\mathcal{D}}$ is open. Indeed, assume to the contrary that there
exist $q_{0}\in\mathcal{A}_{\mathcal{D}}$ and $q_{n}\not \in \mathcal{A}%
_{\mathcal{D}}$ such that $q_{n}\rightarrow q_{0}$. We take nontrivial
solutions $u_{n}\not \gg 0$ of $(P_{\mathcal{D}})$ with $q=q_{n}$. It is
easily seen that $\left\{  u_{n}\right\}  $ is bounded in $H_{0}^{1}(\Omega)$.
Up to a subsequence, $u_{n}\rightarrow u_{0}$ in $C^{1}(\overline{\Omega})$,
where $u_{0}$ is a nontrivial solution of $(P_{\mathcal{D}})$ with $q=q_{0}$.
Since $q_{0}\in\mathcal{A}_{\mathcal{D}}$, we have $u_{0}\gg0$, and so
$u_{n}\gg0$ for $n$ large enough, which is a contradiction. Thus
$\mathcal{A}_{\mathcal{D}}$ is open. The proof of the connectedness of
$\mathcal{A}_{\mathcal{D}}$ is more technical, and we refer to \cite{KRQU16}
for the details.
%\marginpar{added}
The proof of (ii) follows similarly, see also \cite{KRQU16}. \qed \newline

Following a similar strategy, we show that the positivity property also holds
in the Dirichlet case if $a^{-}$ is small enough
%\marginpar{moved to this point and modified}
(assuming now that $q\in(0,1)$ is fixed), which extends Theorem \ref{t11} (ii)
under (A.1). Let us add that this theorem is also true for some non-powerlike
nonlinearities \cite[Theorem 1.1]{KRQU16}.

\begin{theorem}
\label{t3a}

Assume (A.1). Then there exists $\delta>0$ (possibly depending on $q$ and
$a^{+}$) such that every nontrivial nonnegative solution $u$ of
$(P_{\mathcal{D}})$ satisfies that $u\gg0$ if $\left\Vert a^{-}\right\Vert
_{C(\overline{\Omega})}<\delta$.
\end{theorem}

Note that since (A.0) is necessary for the existence of positive solutions of
$(P_{\mathcal{N}})$, we can't expect an analogue of the above theorem for this problem.

As an immediate consequence of Theorem \ref{t2A} and Corollary \ref{A=I}, we infer:

\begin{cor}
\label{old6.2}Assume that $\Omega_{+}$ is connected and satisfies (A.2), and
let $u_{q}$ be the unique nontrivial solution of $(P_{\mathcal{D}})$. Then
$u_{q}\not \gg 0$ for all $q\in\left(  0,q_{\mathcal{D}}\right]  $ and
$u_{q}\gg0$ for all $q\in(q_{\mathcal{D}},1)$. A similar result holds for
$(P_{\mathcal{N}})$ assuming, in addition, (A.0).
\end{cor}

Let us mention that, if in addition to the assumptions of Corollary
\ref{old6.2}, $\Omega_{+}$ includes a \textit{tubular neighborhood of}
$\partial\Omega$ (i.e., a set of the form $\{x\in\Omega:d(x,\partial
\Omega)<\rho\}$, for some $\rho>0$) then the \textbf{SMP} shows that
%\textit{every} nontrivial solution of $(P_{\mathcal{N}})$ is positive on $\partial\Omega$ \cite[Theorem 1.10 (i)]{KRQUnodea}.
the solution $u_{q}$ above satisfies either $u_{q}\gg0$ or $u_{q}=0$ somewhere
in $\Omega$, see Figure \ref{fig04}.\smallskip

Although Theorem \ref{t2A} claims that under (A.1) the sets $\mathcal{A}%
_{\mathcal{D}}$ and $\mathcal{A}_{\mathcal{N}}$ are \textit{always} nonempty,
by Example \ref{exa:cos} we see that given \textit{any} $q\in(0,1)$, we may
find $a=a_{q}$ satisfying (A.1) and such that $(P_{\mathcal{D}})$ and
$(P_{\mathcal{N}})$ have a nontrivial
%\marginpar{removed}
solution $u\not \gg 0$. In view of Theorem \ref{t2A}, this fact shows that
$\mathcal{A}_{\mathcal{D}}$ and $\mathcal{A}_{\mathcal{N}}$ can be arbitrarily
small for a suitable $a$.

The next result (cf. \cite[Theorem 1.4 (i)]{KRQUnodea}, \cite[Proposition 5.1
(i)]{KRQU3}) shows that for any $q\in\left(  0,1\right)  $, we may find $a$
such that $q\in\mathcal{I}(a)\setminus\mathcal{A}(a)$ (and so, in general,
$\mathcal{A}\subsetneq\mathcal{I}$).

\begin{prop}
\strut\label{ti}

\begin{enumerate}
\item Given $\Omega\subset\mathbb{R}$ and $q\in(0,1)$, there exists $a\in
C(\overline{\Omega})$ such that $q\in\mathcal{I}_{\mathcal{N}}\setminus
\mathcal{A}_{\mathcal{N}}$.

\item Given $\Omega\subset\mathbb{R}$ and $q\in(0,1)$, there exists $a\in
C(\Omega)\cap L^{r}\left(  \Omega\right)  $, $r>1$, such that $q\in
\mathcal{I}_{\mathcal{D}}\setminus\mathcal{A}_{\mathcal{D}}$.
\end{enumerate}
\end{prop}

\subsection{The ground state solution}

\label{subsec:ground}

Recall that the Dirichlet eigenvalue problem
%\begin{align}%
\[
\label{Depro}
\begin{cases}
-\Delta\phi=\mu a(x)\phi & \mbox{in}\ \Omega,\\
\phi=0 & \mbox{on}\ \partial\Omega.
\end{cases}
\leqno{(E_{\mathcal{D}})}
\]
has a first positive eigenvalue $\mu_{\mathcal{D}}(a)$, which is principal and
simple, and a positive eigenfunction $\phi_{\mathcal{D}}(a)\gg0$ normalized by
$\int_{\Omega} \phi_{\mathcal{D}}^{2}=1$, associated with $\mu_{\mathcal{D}%
}(a)$.

Let $I_{q}:H_{0}^{1}(\Omega)\rightarrow\mathbb{R}$ be given by
%\marginpar{added $\left(  x\right)  $}%
\[
I_{q}(u):=\frac{1}{2}\int_{\Omega}|\nabla u|^{2}-\frac{1}{q+1}\int_{\Omega
}a(x) |u|^{q+1}%
\]
for $q\in\lbrack0,1)$. It is well-known that nonnegative critical points (in
particular minimizers) of $I_{q}$ are solutions of $(P_{\mathcal{D}})$. 
By a 
%\marginpar{added 'global'}
\textit{ground state} of $I_{q}$ we mean a global
minimizer of this functional.

\begin{prop}
\label{pgs} $I_{q}$ has a unique nonnegative ground state 
%\marginpar{modified}
$U_{q}$ for every $q\in(0,1)$. In addition:

\begin{enumerate}
\item $U_{q}>0$ in $\Omega_{+}$ and $q\mapsto U_{q}$ is continuous from
$(0,1)$ to $W_{\mathcal{D}}^{2,r}(\Omega)$.

\item There exists $q_{0}\in(0,1)$ such that $U_{q}\gg0$ for $q\in(q_{0},1)$.

\item As $q\rightarrow1^{-}$ we have $U_{q}\rightarrow0$ in $C_{0}
^{1}(\overline{\Omega})$ if $\mu_{\mathcal{D}}(a)>1$, whereas $\Vert
U_{q}\Vert_{C(\overline{\Omega})} \rightarrow\infty$ if $\mu_{\mathcal{D}%
}(a)<1$.

\item If $q_{n}\rightarrow0^{+}$ then, up to a subsequence, $U_{q_{n}
}\rightarrow U_{0}$ in $C_{0}^{1}(\overline{\Omega})$, where $U_{0}$ is a
nonnegative global minimizer of $I_{0}$. In particular, if $0\not \equiv
\mathcal{S}(a)\geq0$ in $\Omega$, then $U_{q}\rightarrow\mathcal{S}(a)$ in
$C_{0}^{1}(\overline{\Omega})$ as $q\rightarrow0^{+}$.
\end{enumerate}
\end{prop}

\textit{Sketch of the proof}. By a standard minimization argument, one may
easily prove the existence of a global minimizer of $I_{q}$. Moreover, there
is a $1$ to $1$ correspondence between global minimizers of $I_{q}$ and
minimizers of $\int_{\Omega}|\nabla u|^{2}$ over the $C^{1}$ manifold
$\left\{  u \in H_{0}^{1}(\Omega): \int_{\Omega}a(x)|u|^{q+1}=1\right\}  $. By
\cite[Theorem 1.1]{KLP}, we infer that if $U_{q}$ and $V_{q}$ are global
minimizers of $I_{q}$ then $U_{q}=tV_{q}$ for some $t>0$. But since $U_{q}$
and $V_{q}$ solve $(P_{\mathcal{D}})$, we deduce that $t=1$, i.e. $U_{q}$ is
the unique nonnegative global minimizer of $I_{q}$. If $U_{q}(x)=0$ for some
$x\in\Omega_{+}$ then, by the \textbf{SMP}, $U_{q}$ vanishes is some ball
$B\subset\Omega_{+}$. We choose a nontrivial and smooth $\psi\geq0$ supported
in $B$ and extend it by zero to $\Omega$. Then $I_{q}(U_{q}+t\psi)=I_{q}%
(U_{q})+I_{q}(t\psi)< I_{q}(U_{q})$ if $t$ is small enough, which yields a
contradiction. Using standard compactness arguments and the uniqueness of
$U_{q}$, we can show that $U_{q} \to U_{q_{0}}$ in $W_{\mathcal{D}}%
^{2,r}(\Omega)$ as $q \to q_{0}$, for any $q_{0}\in(0,1)$. Arguing as in the
proof of Theorem \ref{t2A} we prove that $U_{q}\gg0$ for $q$ close to $1$, and
$U_{q}\rightarrow0$ in $C_{0}^{1}(\overline{\Omega})$ if $\mu_{\mathcal{D}%
}(a)>1$. If $\mu_{\mathcal{D}}(a)<1$ and $\{u_{n}\}$ is bounded in $H_{0}%
^{1}(\Omega)$, where $u_{n}:=U_{q_{n}}$ and $q_{n} \to1^{-}$, then again as in
the proof of Theorem \ref{t2A}, we find that $u_{n} \to u_{0}$ and $u_{0}
\geq0$ solves $-\Delta u_{0}=a(x)u_{0}$ in $\Omega$, $u_{0}=0$ on
$\partial\Omega$. Using the fact that $u_{n}$ are ground state solutions, we
can show that $u_{0} \not \equiv 0$, so that $\mu_{\mathcal{D}}(a)=1$, a
contradiction. Finally, we refer to \cite{KRQU3} for the proof of (iv).
\qed \newline

\begin{rem}
\strut

\begin{enumerate}
\item It is not hard to show that under (A.0) the functional $I_{q}$,
considered now in $H^{1}(\Omega)$, has a ground state, which is positive in
$\Omega_{+}$, and strongly positive for $q$ close enough to $1$.

\item Proposition \ref{pgs} (ii) extends Theorem \ref{t11} and Theorem
\ref{t2A}(i) (as long as the existence of a positive solution is concerned)
without assuming (A.1).
%\marginpar{added}

\end{enumerate}
\end{rem}

\section{Structure of the positive solutions set}

%\section{Bifurcation for Dirichlet and Neumann conditions}

\label{sec:BDN}

This section is devoted to a further investigation of the set $\mathcal{I}%
_{\mathcal{D}}$ (respect.\ $\mathcal{I}_{\mathcal{N}}$), which provides a
rather complete description of
%\marginpar{paragraph modified}
the positive solutions set of $(P_{\mathcal{D}})$ (respect.\ $\left(
P_{\mathcal{N}}\right)  $). From Theorem \ref{t2A} we observe that
$(q_{\mathcal{D}},1)\subseteq\mathcal{I}_{\mathcal{D}}$ and $(q_{\mathcal{N}%
},1)\subseteq\mathcal{I}_{\mathcal{N}}$. Taking advantage of the ground state
solution, constructing suitable sub and supersolutions, and also using a
bifurcation approach, we analyze the asymptotic behavior of the positive
solutions as $q\rightarrow1^{-}$ and $q\rightarrow0^{+}$.

%\marginpar{\textit{{\footnotesize removed}}}
%We remark that our results in this section remain valid for the more general situation that $a\in L^{r}(\Omega)$ with $r>N$, where we assume (A.1$^{+}$) instead of (A.1).

\subsection{The Dirichlet problem}

\label{subsec:D}

Let us consider %\marginpar{modified}
$(P_{\mathcal{D}})$, with $q\in(0,1)$ as
a bifurcation parameter. To this end, we introduce two further conditions on
$a$. The first one slightly weakens (A.3) requiring that
\[
\mathcal{S}(a)>0\ \ \mbox{in}\ \ \Omega,\leqno{({\bf A.3^{\prime}})}
\]
whereas the second one is a technical decay condition near $\partial\Omega$:
%\[
%\left\{
%\begin{array}
%[c]{l}%
%\Omega_{+}\text{ is connected and }\\
%\partial\Omega_{+}\text{ satisfies an inner sphere condition with respect to
%}\Omega_{+},
%\end{array}
%\right.  \leqno{(H_{+}')}
%\]%
\[
\left\vert a(x)\right\vert \leq Cd(x,\partial\Omega)^{\eta}\text{
\ \negthinspace a.e.\ in }\Omega_{\rho_{0}},\text{ for some }\rho_{0}>0\text{
and }\eta>1-\frac{1}{N},\leqno{({\bf A.4})}
\]
where
\begin{equation}
\Omega_{\rho}=:\{x\in\Omega:d(x,\partial\Omega)<\rho\}. \label{def:tnei}%
\end{equation}
Recall
%\marginpar{\textit{added}}
that $\Omega_{\rho}$ is said to be a \textit{tubular neighborhood of}
$\partial\Omega$. It turns out that (A.3$^{\prime}$) is sufficient to deduce
the conclusion of Theorem \ref{t11} (i), i.e. that $(P_{\mathcal{D}})$ has a
positive solution for every $q\in(0,1)$. In addition, we shall use
(A.3$^{\prime}$) to show that this solution converges to $\mathcal{S}(a)$ as
$q\rightarrow0^{+}$. On the other hand, (A.4) is needed to obtain solutions of
$(P_{\mathcal{D}})$ bifurcating from $t\phi_{\mathcal{D}}$, for some $t>0$,
when $\mu_{\mathcal{D}}(a)=1$. Since
%\marginpar{modified}
$\phi_{\mathcal{D}}=0$ on $\partial\Omega$, we assume (A.4) to ensure that
$a\,\phi_{\mathcal{D}}^{q-2}$ has the appropriate integrability to carry out
this bifurcation procedure, see Subsection \ref{ssubsec:IFTLB}.

\begin{figure}[tbh]
\centerline{
		\includegraphics[scale=0.12]{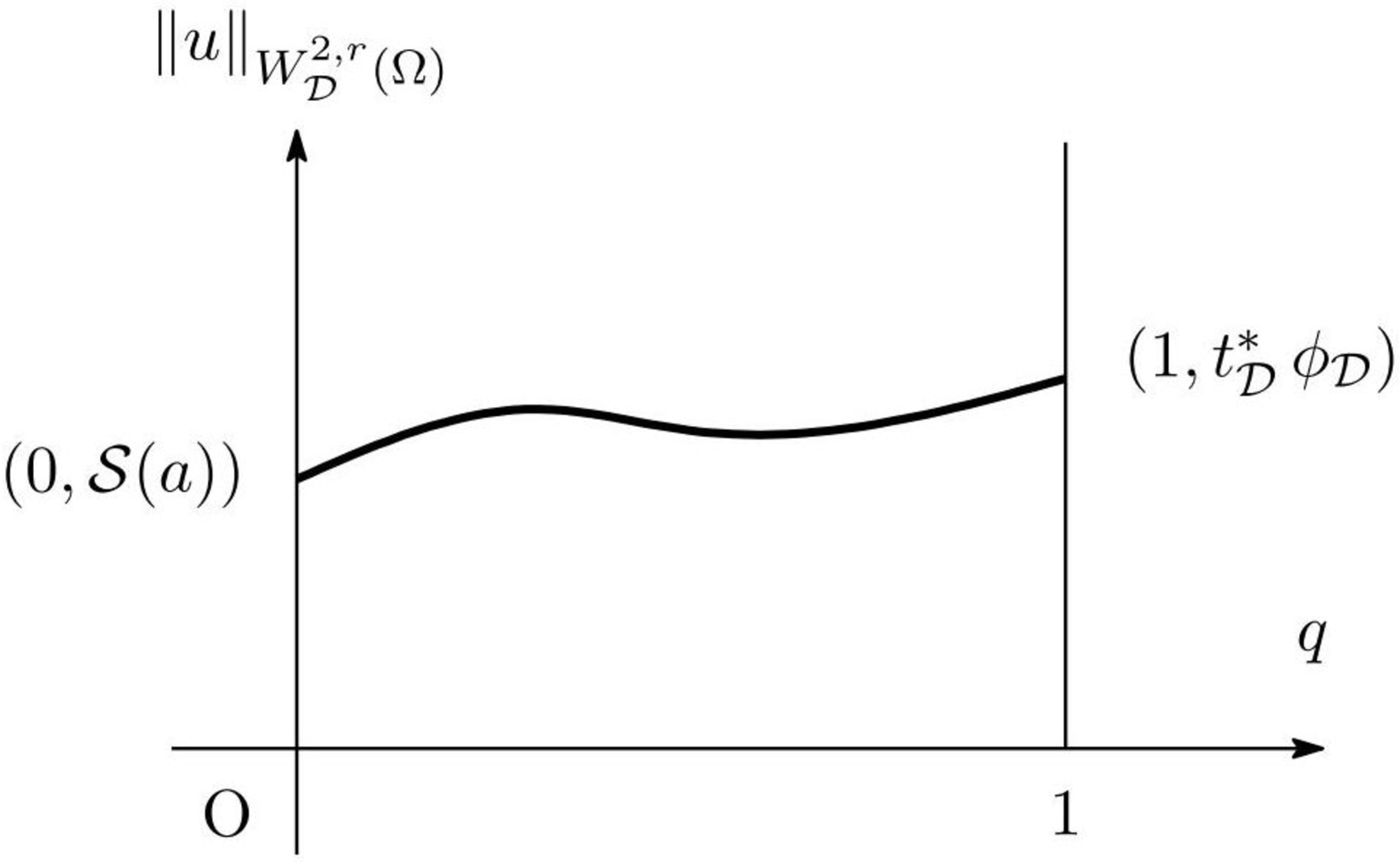}
		\hskip0.2cm
		\includegraphics[scale=0.12]{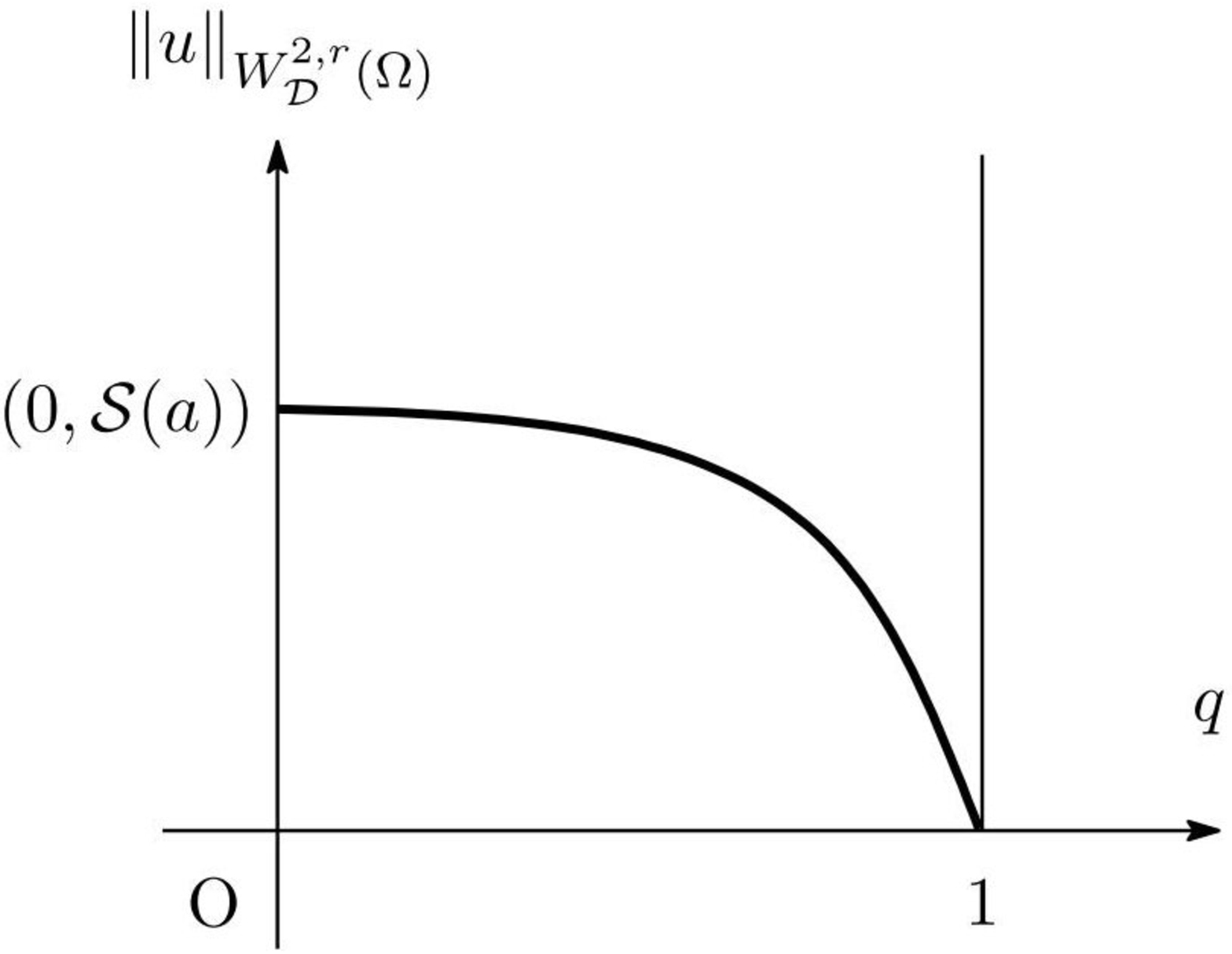}
		\hskip0.0cm
		\includegraphics[scale=0.12]{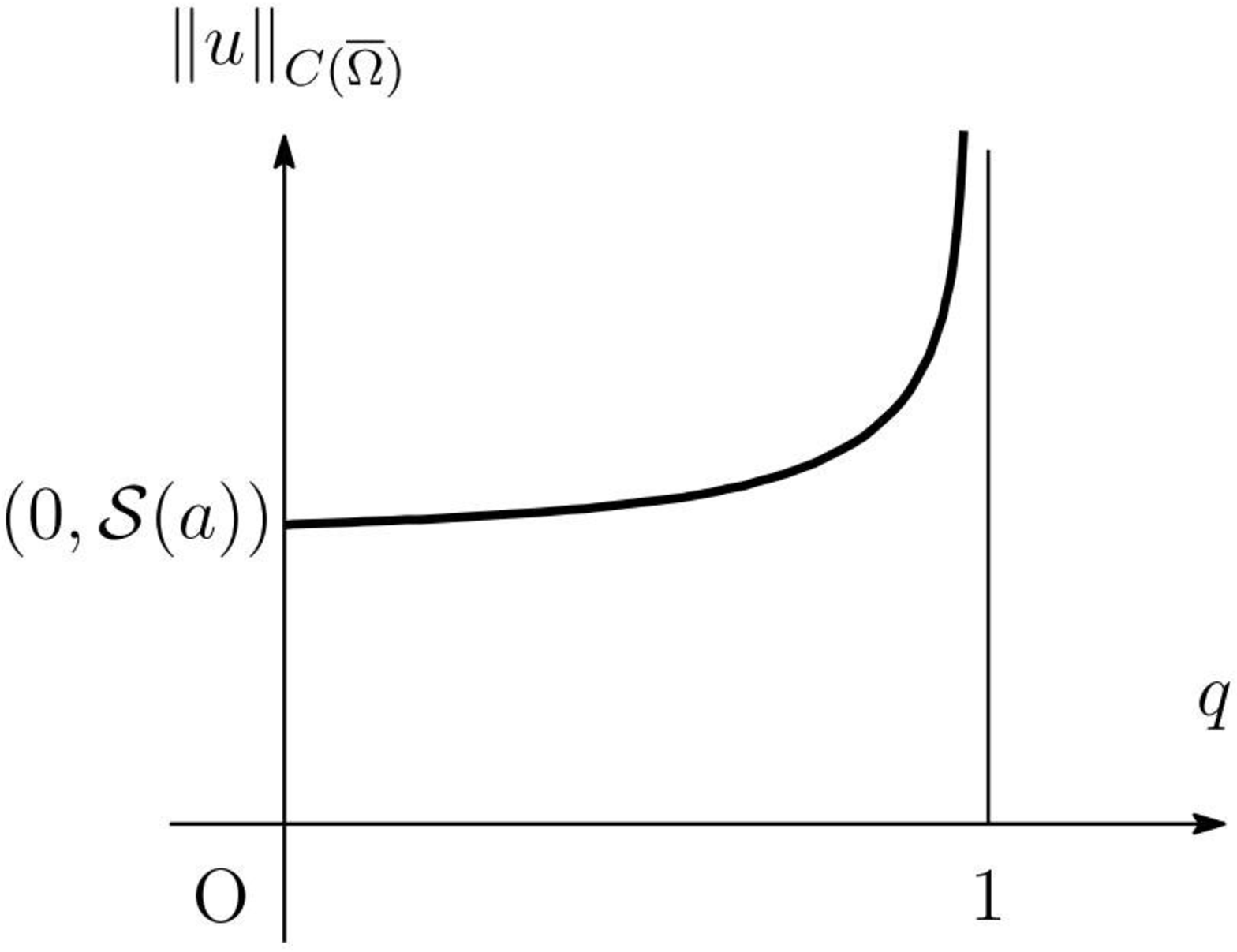}
	} \centerline{(i) \hskip3.7cm (ii) \hskip3.7cm (iii)}\caption{The curve of
positive solutions emanating from $(0,\mathcal{S}(a))$: Cases (i)
$\mu_{\mathcal{D}}(a)=1$, (ii) $\mu_{\mathcal{D}}(a)>1$, (iii) $\mu
_{\mathcal{D}}(a)<1$.}%
\label{fig:double}%
\end{figure}

Denoting by $u_{\mathcal{D}}(q)$ the unique positive solution of
$(P_{\mathcal{D}})$ for $q\in(0,1)$ whenever it exists, we see from
%\marginpar{fixed typo twice;}
Proposition \ref{pgs} (ii) that $U_{q}=u_{\mathcal{D}}(q)$ for $q$ close to
$1$, so that Proposition \ref{pgs} (iii) provides the asymptotics of
$u_{\mathcal{D}}(q)$ when $\mu_{\mathcal{D}}(a)\neq1$.
%\marginpar{modified this part and added sentence}
We treat now the case $\mu_{\mathcal{D}}(a)=1$ and also provide the asymptotic
behavior of $u_{\mathcal{D}}(q)$ as $q\rightarrow0^{+}$, as well as sufficient
conditions to have $u_{\mathcal{D}}(q)\gg0$ for \textit{every} $q\in(0,1)$.
%\marginpar{\textit{removed}}
Under these conditions, we obtain a rather complete description of the
positive solutions set $\left\{  (q,u_{\mathcal{D}}(q)):q\in(0,1)\right\}  $
of $(P_{\mathcal{D}})$, see Figure \ref{fig:double}. We shall present here a
simplified version of these results. For the precise assumptions required in
each of following items we refer to \cite[Theorems 1.2 and 1.4, Corollary
1.6]{KRQU3}. Under (A.4), let us set
\begin{equation}
t_{\mathcal{D}}^{\ast}:=\exp\left[  -\frac{\int_{\Omega}a(x)\phi_{\mathcal{D}%
}^{2}\log\phi_{\mathcal{D}}}{\int_{\Omega}a(x)\phi_{\mathcal{D}}^{2}}\right]
. \label{def:tdast}%
\end{equation}

\begin{theorem}
\label{cormt2}

Let $r>N$. Assume (A.1), (A.2), (A.3$\,^{\prime}$) and (A.4). Then
$u_{\mathcal{D}}(q)=U_{q}>0$ in $\Omega$ for every $q\in(0,1)$. In addition,
if we set $u_{\mathcal{D}}(0):=\mathcal{S}(a)$ then $q\mapsto u_{\mathcal{D}%
}(q)$ is continuous from $[0,1)$ to $W_{\mathcal{D}}^{2,r}(\Omega)$. The
asymptotic behavior of $u_{\mathcal{D}}(q)$ as $q\rightarrow1^{-}$ is
characterized as follows:

\begin{enumerate}
\setlength{\itemsep}{0.0cm}

\item If $\mu_{\mathcal{D}}(a)\geq1$ and we set
\[
u_{\mathcal{D}}(1):= \left\{
\begin{array}
[c]{ll}%
t_{\mathcal{D}}^{\ast}\,\phi_{\mathcal{D}}, & \mbox{if} \ \ \mu_{\mathcal{D}
}(a)=1,\\
0, & \mbox{if} \ \ \mu_{\mathcal{D}}(a)>1 \ (\mbox{bifurcation from zero}),
\end{array}
\right.
\]
then $q\mapsto u_{\mathcal{D}}(q)$ is left continuous at $q=1$ (see Figure
\ref{fig:double} (i), (ii)).
%Moreover, if in addition some of the conditions in Theorem \ref{mt2b} (iii) hold, then $u_{D}(q)\gg0$ for all $q\in(0,1)$.

\item If $\mu_{\mathcal{D}}(a)<1$ then the curve $\{ (q,u_{\mathcal{D}}(q)): q
\in[0,1) \}$
%of positive solutions
bifurcates from infinity at $q=1$ (see Figure \ref{fig:double} (iii)).
\end{enumerate}

Finally, as for the strong positivity of $u_{\mathcal{D}}(q)$, we have the
following two assertions:

\begin{enumerate}
\setlength{\itemsep}{0.0cm}

\item[(iii)] If (A.3) holds then $u_{\mathcal{D}}(q)\gg0$ for $q$ close to $0$
or $1$.

\item[(iv)] In the following cases, we have $u_{\mathcal{D}}(q)\gg0$ for
\textrm{all} $q\in(0,1)$
%\marginpar{added}
(and so, $\mathcal{I}_{\mathcal{D}}=(0,1)$):

\begin{enumerate}
\setlength{\itemsep}{0.2cm}

\item $a\geq0$ in $\Omega_{\rho_{0}}$ for some $\rho_{0}>0$,

\item $\Omega$ is a ball and $a$ is radial,

\item (A.3) holds
%$\mathcal{S}(a)\gg0$,
and $\Omega_{+}$ is connected.
%\marginpar{\textit{removed}}
%\newline

\end{enumerate}
\end{enumerate}
\end{theorem}

%\begin{figure}[tbh]
%\label{f1}
%\par
%\begin{center}
%\includegraphics[scale=0.22]{fig17_0618b-m16nov19.eps}
%\end{center}
%\caption{The curve of positive solutions in the case %$\mu_{\mathcal{D}}
%(a)=1$.}%
%\label{fig17_0618b}%
%\end{figure}

\begin{rem}
\strut\label{remgs}

\begin{enumerate}
\item From Theorem \ref{t10} and
%\marginpar{fixed typo}
Proposition \ref{pgs} (ii), it suffices to assume (A.1) and (A.2) to have
$U_{q}=u_{\mathcal{D}}(q)$ whenever $u_{\mathcal{D}}(q)$ exists. Moreover,
under these conditions,
\begin{equation}
\mathcal{I}_{\mathcal{D}}=\{q\in(0,1):U_{q}\gg0\}, \label{iia}%
\end{equation}
and $\mathcal{I}_{\mathcal{D}}$ is open.

%\item Whenever $U_{q}>0$ in $\Omega$, we have, by Theorem \ref{t10} (i-a), $U_{q}=u_{\mathcal{D}}(q)$. This equality also holds whenever $u_{\mathcal{D}}(q)$ exists and (A.1) and (A.2) are satisfied, as a consequence of the first assertion in Theorem \ref{mt2a} (iii).

\item The assertion in Theorem \ref{cormt2} (i) when $\mu_{\mathcal{D}}(a)=1$
also gives a better asymptotic estimate for
%\marginpar{modified}
$U_{q}$ as $q\rightarrow1^{-}$ if (A.4) holds and $\mu_{\mathcal{D}}(a)\neq1$.
Indeed, a rescaling argument yields that
\[
U_{q}\sim\mu_{1}(a)^{-\frac{1}{1-q}}\,t_{\mathcal{D}}^{\ast}\,\phi
_{\mathcal{D}}\quad\mbox{as}\quad q\rightarrow1^{-},
\]
i.e.
\[
\mu_{\mathcal{D}}(a)^{\frac{1}{1-q}}U_{q}\rightarrow t_{\mathcal{D}}^{\ast
}\,\phi_{\mathcal{D}}\ \mbox{ in }\ W_{\mathcal{D}}^{2,r}(\Omega
)\ \mbox{ as }\ q\rightarrow1^{-}.
\]

\item As already stated, under (A.3$^{\prime}$) we have a positive solution
%\marginpar{removed}
for every $q\in(0,1)$. Assuming additionally (A.3), we can deduce the
conclusion of Theorem \ref{cormt2} (iii), which extends Theorem \ref{t11}
(i).
%\marginpar{added}
Let us add that in some cases, by Theorem \ref{cc} below, the condition
$\int_{\Omega}a\geq0$ (which is weaker than (A.3$^{\prime}$)) is also
sufficient to have a positive solution of $\left(  P_{\mathcal{D}}\right)  $
for all $q\in(0,1)$.

\item Under the assumptions of Theorem \ref{cormt2} (iv-c), we infer from
Corollary \ref{A=I} that $\mathcal{A}_{\mathcal{D}}=(0,1)$.
%, see Corollary \ref{cor:0to1}.

\end{enumerate}
\end{rem}

%As an immediate consequence of Theorem \ref{mt2b} (iii-c) and , we have the following result.

%\begin{cor}
%\label{A1}
%\end{cor}

Next we consider the linearized stability of a solution in $\mathcal{P}%
_{\mathcal{D}}^{\circ}$ of $(P_{\mathcal{D}})$ for $q\in\mathcal{I}%
_{\mathcal{D}}$. Let us recall that a solution $u\gg0$ of $(P_{\mathcal{D}})$
is said to be \textit{asymptotically stable} if $\gamma_{1}(q,u)>0$, where
$\gamma_{1}(q,u)$ is the first eigenvalue of the linearized eigenvalue problem
at $u$, namely,%
\begin{equation}%
\begin{cases}
-\Delta\varphi=qa(x)u^{q-1}\varphi+\gamma\varphi & \mbox{ in }\Omega,\\
\varphi=0 & \mbox{ on }\partial\Omega.
\end{cases}
\label{stab}%
\end{equation}
Observe that under the decay condition (A.4), given $q\in\left[  0,1\right)  $
and $u\gg0$, we have $au^{q-1}\in L^{t}\left(  \Omega\right)  $ for some
$t>N$, so that $\gamma_{1}(q,u)$ is well defined.

The implicit function theorem (\textbf{IFT} for short) provides us with the
following result \cite[Theorem 1.5]{KRQU3}:

\begin{theorem}
\label{mt2c} If (A.4) holds then $\mathcal{I}_{\mathcal{D}}$ is open, and
$u_{\mathcal{D}}(q)$ is asymptotically stable for $q\in\mathcal{I}%
_{\mathcal{D}}$.
\end{theorem}

\subsubsection{Local bifurcation analysis in the case $\mu_{\mathcal{D}}%
(a)=1$}

\label{ssubsec:IFTLB}

Let us give a sketch of the proof of Theorem \ref{cormt2} (i) when
$\mu_{\mathcal{D}}(a)=1$. In this case, $(P_{\mathcal{D}})$ has the trivial
line of strongly positive solutions:
\[
\Gamma_{1}:=\left\{  (q,u)=(1,t\phi_{\mathcal{D}}):t>0\right\}  .
\]
For $q\simeq1$, where $q$ is a bifurcation parameter, we shall construct
solutions of $(P_{\mathcal{D}})$ bifurcating at certain $(1,t\phi
_{\mathcal{D}})\in\Gamma_{1}$ in $\mathbb{R}\times W_{\mathcal{D}}^{2,\xi
}(\Omega)$, for some fixed $\xi>N$. This bifurcation result (Proposition
\ref{prop:gam2} below) complements Proposition \ref{pgs} (iii).

Under (A.4), choose $\sigma_{0}>0$ such that $\eta>1+\sigma_{0}-\frac{1}{N}$
and set $J_{0}:=(1-\frac{\sigma_{0}}{2},1+\frac{\sigma_{0}}{2})$. We fix then
$\xi\in(N,r)$, depending only on $N$ and $\sigma_{0}$, in such a way that
$\xi(\eta+q-2)>-1+\frac{\sigma_{0}N}{4}$ for $q\in J_{0}$. Following the
Lyapunov-Schmidt procedure, we reduce $(P_{\mathcal{D}})$ to a bifurcation
equation. Set $A:=-\Delta-a(x)$ with domain $D(A):=W_{\mathcal{D}}^{2,\xi
}(\Omega)$. Then $\mathrm{Ker}A=\left\{  t\phi_{\mathcal{D}} :t\in\mathbb{R}
\right\}  $ and $\mathrm{Im}A= \left\{  f \in L^{\xi}(\Omega):\int_{\Omega}f
\phi_{\mathcal{D}}=0 \right\}  $. Let $Q$ be the projection of $L^{\xi}%
(\Omega)$ to $\mathrm{Im}A$, given by $Q[f]:=f-\left(  \int_{\Omega}%
f\phi_{\mathcal{D}} \right)  \phi_{\mathcal{D}}$. As long as we consider
solutions $u\gg0$, $(P_{\mathcal{D}})$ is equivalent to the following coupled
equations: for $u=t\phi_{\mathcal{D}}+w\in D(A)=\mathrm{Ker}A+X_{2}$ with
$t=\int_{\Omega}u\phi_{\mathcal{D}}$ and $X_{2}=\{ u\in D(A):\int_{\Omega
}u\phi_{\mathcal{D}}=0 \}$,
\begin{align}
&  Q\left[  A(t\phi_{\mathcal{D}}+w) \right]  = Q\left[  a\left(  x\right)
\left(  (t\phi_{\mathcal{D}}+w)^{q}-(t\phi_{\mathcal{D}}+w) \right)  \right]
,\label{eq:Q}\\
&  (1-Q)\left[  A(t\phi_{\mathcal{D}}+w) \right]  =(1-Q)\left[  a\left(
x\right)  \left(  (t\phi_{\mathcal{D}}+w)^{q}-(t\phi_{\mathcal{D}}+w) \right)
\right]  . \label{eq:1-Q}%
\end{align}

Given $t_{0}>0$, first we solve \eqref{eq:Q} with respect to $w$ at
$(q,t,w)=(1,t_{0},0)$, where $(1,t_{0},0)$ is a solution of \eqref{eq:Q}. Note
that (A.4) gives that \eqref{eq:Q} is $C^{2}$ for $(q,t,w)\simeq(1,t_{0},0)$,
since the choice of $\xi$ ensures that $a(t\phi_{\mathcal{D}}+w)^{q-2}\in
L^{\xi}(\Omega)$ for such $(q,t,w)$. An \textbf{IFT} argument shows the
existence of a unique $w=w(q,t)$ for every $(q,t)\simeq(1,t_{0})$ such that
$(q,t,w)$ solves \eqref{eq:Q}. We plug $w(q,t)$ into \eqref{eq:1-Q}, and thus,
deduce the desired bifurcation equation
\[
\Phi(q,t):=\int_{\Omega}a(x)\{(t\phi_{\mathcal{D}}+w(q,t))^{q} -(t\phi
_{\mathcal{D}}+w(q,t)\}\phi_{\mathcal{D}}=0,\quad(q,t) \simeq(1,t_{0}),
\]
where we note that $\Phi$ is $C^{2}$ for $(q,t)\simeq(1,t_{0})$.

As an application of the \textbf{IFT}, we find that if $(1,t_{0}%
\phi_{\mathcal{D}})$ is a bifurcation point on $\Gamma_{1}$ then
\[
\frac{\partial\Phi}{\partial q}(1,t_{0})=t_{0}\left\{  (\log t_{0}%
)\int_{\Omega}a(x)\phi_{\mathcal{D}}^{2}+\int_{\Omega}a(x)\phi_{\mathcal{D}%
}^{2}\log\phi_{\mathcal{D}}\right\}  =0,
\]
so that $t_{0}=t_{\mathcal{D}}^{\ast}$, given by \eqref{def:tdast}.
Conversely, since direct computations \cite[Lemma 4.3]{KRQU3} provide
\[
\frac{\partial\Phi}{\partial t}(1,t_{\mathcal{D}}^{\ast})=\frac{\partial
^{2}\Phi}{\partial t^{2}}(1,t_{\mathcal{D}}^{\ast})=0,\quad\frac{\partial
^{2}\Phi}{\partial t\partial q}(1,t_{\mathcal{D}}^{\ast})=\int_{\Omega
}a\left(  x\right)  \phi_{\mathcal{D}}^{2}>0,
\]
the Morse Lemma \cite[Theorem 4.3.19]{DM13} yields the following existence
result \cite[Proposition 4.4]{KRQU3}:

\begin{prop}
\label{prop:gam2}

Suppose (A.4) with $\mu_{\mathcal{D}}(a)=1$. Then the set of solutions of
$(P_{\mathcal{D}})$ near $(1,t_{\mathcal{D}}^{\ast}\phi_{\mathcal{D}})$
consists of two continuous curves in $\mathbb{R}\times W_{\mathcal{D}}^{2,\xi
}(\Omega)$ intersecting only at $(1,t_{\mathcal{D}}^{\ast}\phi_{\mathcal{D}})$
\textrm{transversally}, given by $\Gamma_{1}\cup\Gamma_{2}$, where $\Gamma
_{2}$ for $q<1$ represents the ground state solution $U_{q}$.
\end{prop}

Let us mention that Proposition \ref{prop:gam2} remains true in $\mathbb{R}%
\times W_{\mathcal{D}}^{2,r}(\Omega)$ by elliptic regularity.

\subsection{The Neumann problem}

\label{subsec:N}
%\marginpar{\textit{maybe we could rename the section, since there are many results concerning }$\left(  P_{\mathcal{D}}\right)  $}

Under (A.0), the Neumann eigenvalue problem
\[%
\begin{cases}
-\Delta\phi=\mu a(x)\phi & \mbox{in}\ \Omega,\\
\partial_{\nu}\phi=0 & \mbox{on}\ \partial\Omega.
\end{cases}
\leqno{(E_{\mathcal{N}})}
\]
has a first positive eigenvalue $\mu_{\mathcal{N}}(a)$, which is principal and
simple, and an eigenfunction $\phi_{\mathcal{N}}(a)\gg0$ associated to
$\mu_{\mathcal{N}}(a)$ and satisfying $\int_{\Omega}\phi_{\mathcal{N}}^{2}=1$.

\begin{figure}[tbh]
\centerline{
\includegraphics[scale=0.15]{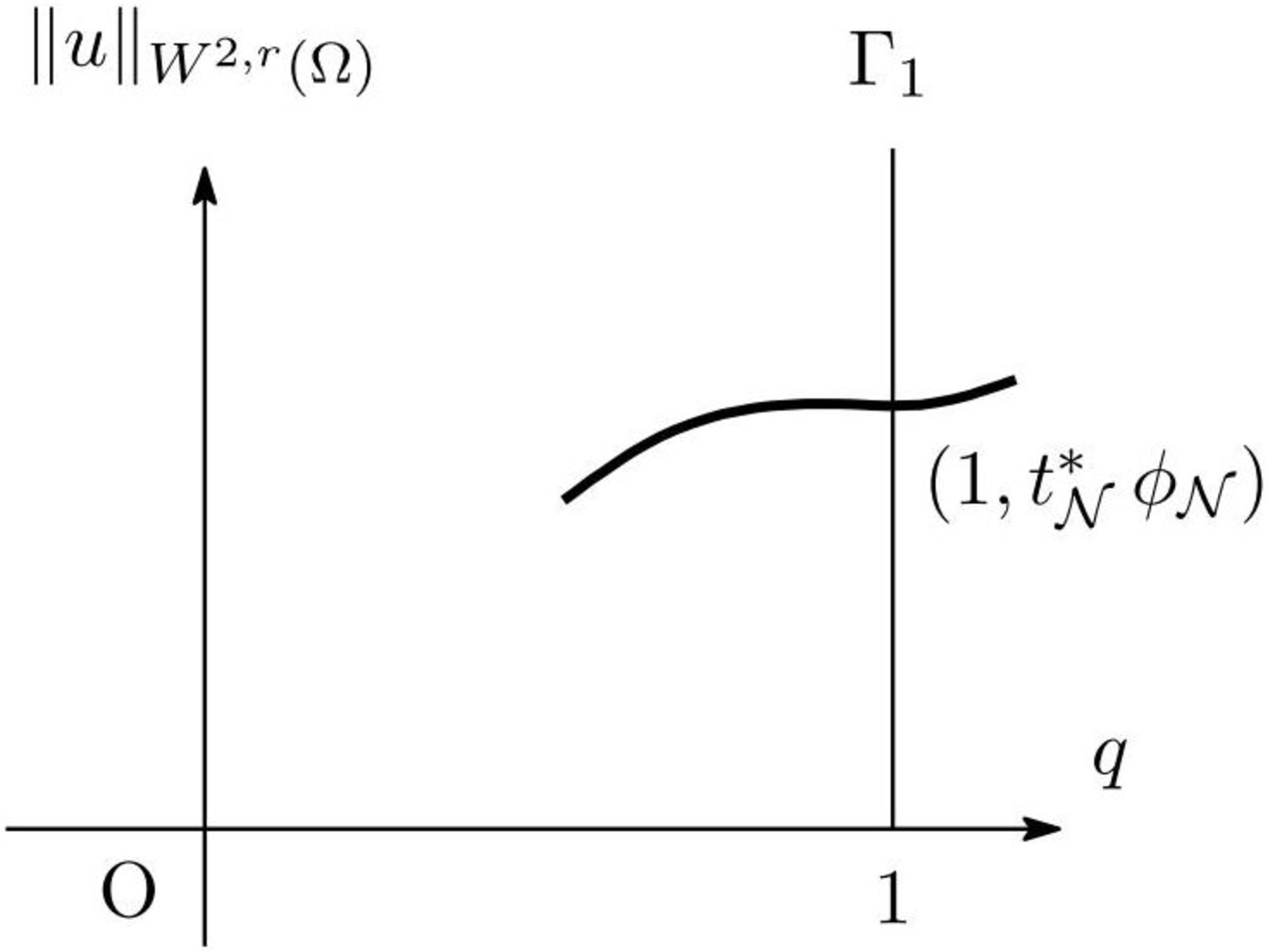}
\hskip0.15cm
\includegraphics[scale=0.15] {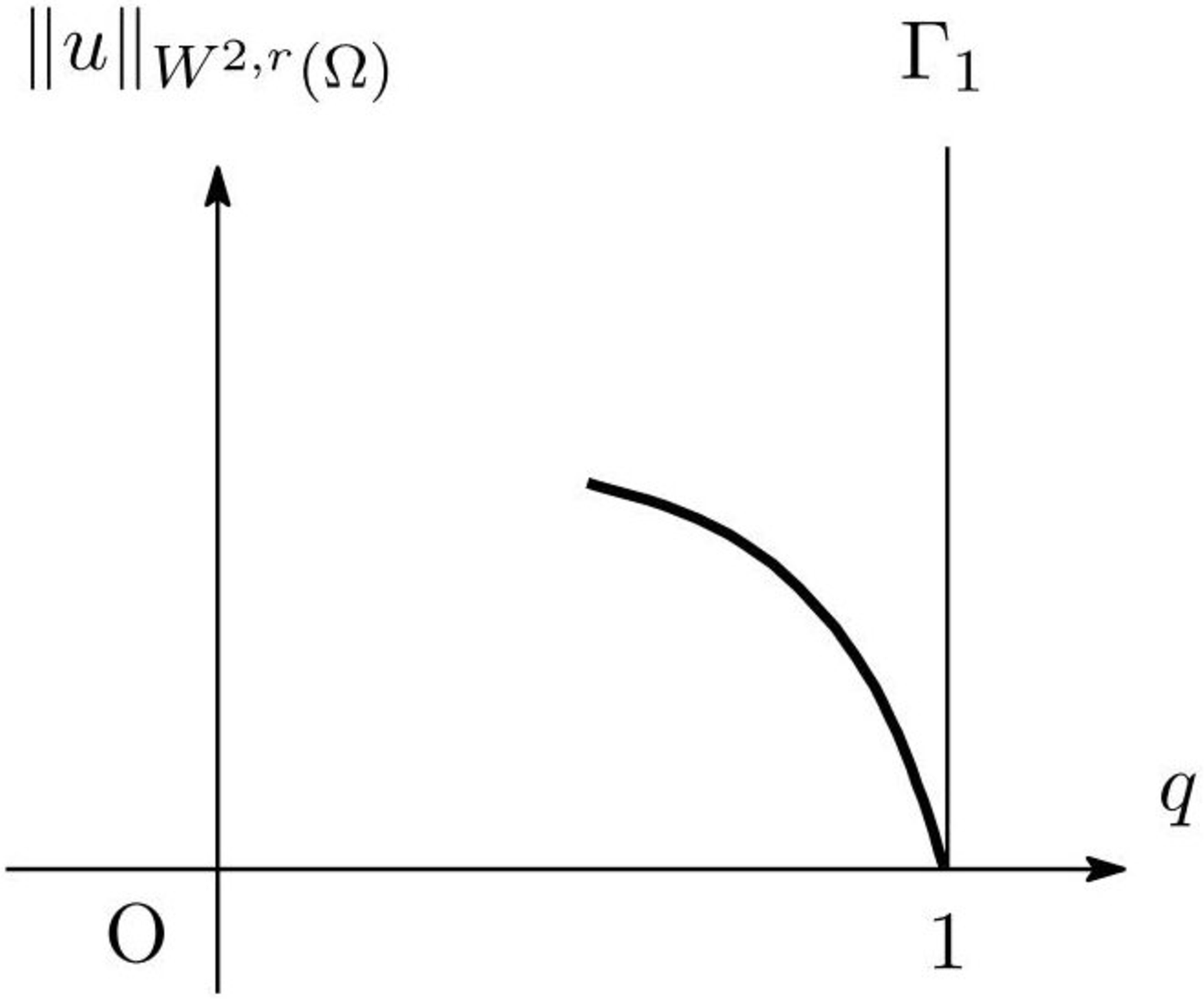}
\hskip0.15cm
\includegraphics[scale=0.15] {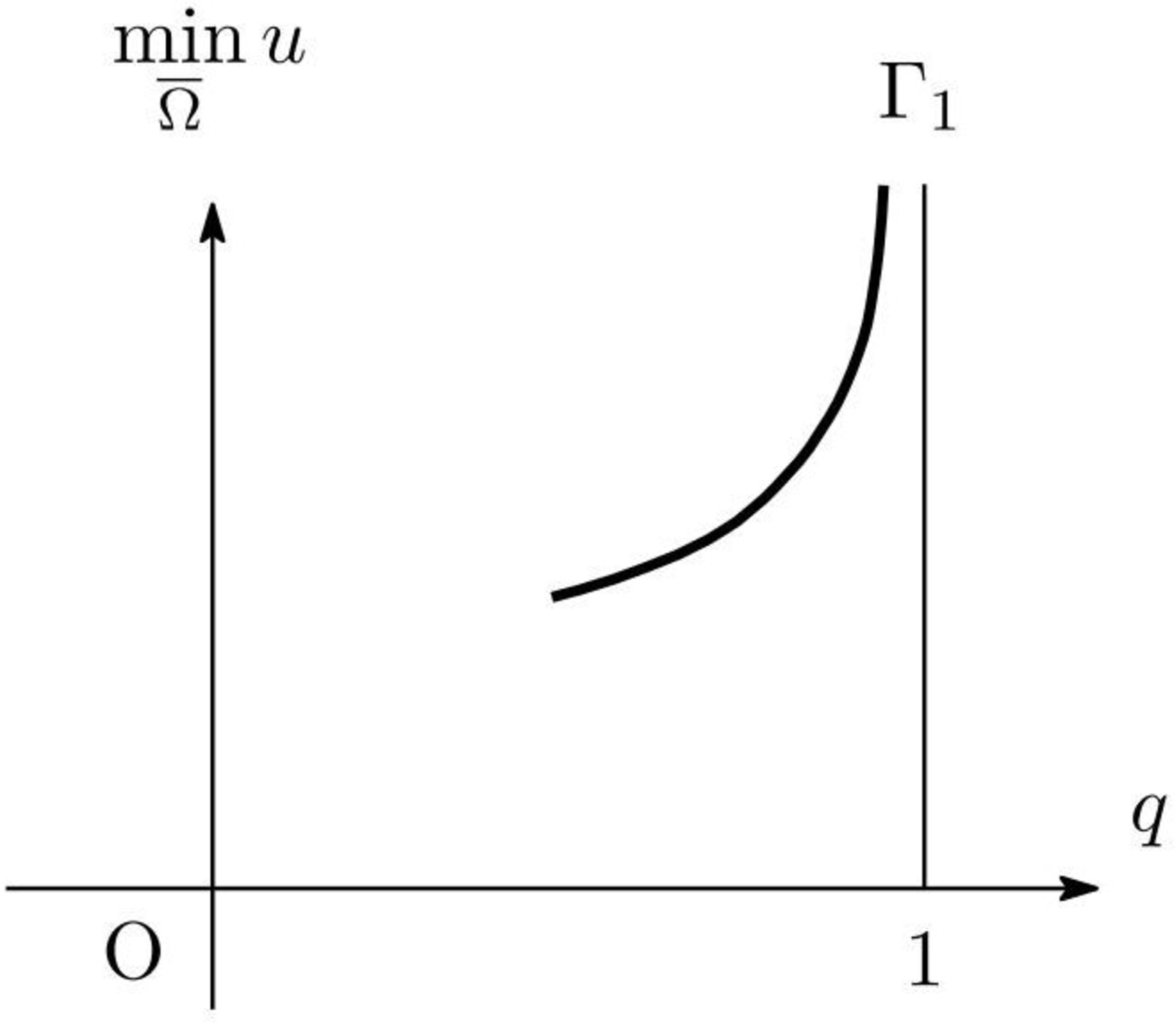}
} \centerline{(i) \hskip3.5cm (ii) \hskip3.5cm (iii)}\caption{Bifurcating
solutions $u\gg0$ (i) from $\Gamma_{1}$ at $\left(  1, t_{\mathcal{N}}^{\ast
}\, \phi_{\mathcal{N}} \right)  $ in case $\mu_{\mathcal{N}}(a)=1$; (ii) from
zero in case $\mu_{\mathcal{N}}(a)>1$; (iii) from infinity in case
$\mu_{\mathcal{N}}(a)<1$.}%
\label{figbif}%
\end{figure}

The bifurcation scheme from the previous subsection also applies to
$(P_{\mathcal{N}})$, with the advantage of \textit{not }requiring the decay
condition (A.4), since %\marginpar{modified}
$\phi_{\mathcal{N}}>0$ on 
$\overline{\Omega}$. We look at $q$ as a bifurcation parameter in
$(P_{\mathcal{N}})$. Similarly as in the Dirichlet case, if $\mu_{\mathcal{N}%
}(a)=1$ then $u=t\phi_{\mathcal{N}}$ solves $(P_{\mathcal{N}})$ with $q=1$,
i.e., $(P_{\mathcal{N}})$ has the trivial line
\[
\Gamma_{1}:=\{(q,u)=(1,t\phi_{\mathcal{N}}):t>0\}.
\]
We shall obtain, for $q$ close to $1$, a curve of solutions $u\gg0$
bifurcating from $\Gamma_{1}$ (see Figure \ref{figbif}).

The definition of \textit{asymptotically stable }for solutions $u\gg0$ of
$(P_{\mathcal{N}})$ is similar to the one for $\left(  P_{\mathcal{D}}\right)
$, see (\ref{stab}). Setting
\begin{equation}
t_{\mathcal{N}}^{\ast}:=\exp\left[  -\frac{\int_{\Omega}a(x)\phi_{\mathcal{N}%
}^{2}\log\phi_{\mathcal{N}}}{\int_{\Omega}a(x)\phi_{\mathcal{N}}^{2}}\right]
, \label{et*}%
\end{equation}
we have the following result \cite[Theorem 1.2]{KRQUnodea}.

\begin{theorem}
\label{c1} Assume (A.0) and $r>N$. Then there exists $q_{0}=q_{0}(a)\in(0,1)$
such that $(P_{\mathcal{N}})$ has a solution $u_{q}\gg0$ for $q_{0}<q<1$.
Moreover, $u_{q}$ is asymptotically stable and satisfies
\[
u_{q}\sim\mu_{\mathcal{N}}(a)^{-\frac{1}{1-q}}\,t_{\mathcal{N}}^{\ast}%
\,\phi_{\mathcal{N}}\quad\mbox{as}\quad q\rightarrow1^{-},
\]
i.e. $\mu_{\mathcal{N}}(a)^{\frac{1}{1-q}}u_{q}\rightarrow t_{\mathcal{N}%
}^{\ast}\,\phi_{\mathcal{N}}$ in $W^{2,r}(\Omega)$ as $q\rightarrow1^{-}$.
More specifically (see Figure \ref{figbif}):

\begin{enumerate}
\setlength{\itemsep}{0.2cm}

\item If $\mu_{\mathcal{N}}(a)=1$, then $u_{q}\rightarrow t_{\mathcal{N}%
}^{\ast}\, \phi_{\mathcal{N}}$ in $W^{2,r}(\Omega)$ as $q\rightarrow1^{-}$.

\item If $\mu_{\mathcal{N}}(a)>1$, then $u_{q}\rightarrow0$ in $W^{2,r}%
(\Omega)$ as $q\rightarrow1^{-}$.

\item If $\mu_{\mathcal{N}}(a)<1$, then $\displaystyle\min_{\overline{\Omega}%
}u_{q}\rightarrow\infty$ as $q\rightarrow1^{-}$.
\end{enumerate}
\end{theorem}

Let us point out that, in general, it is hard to give a lower estimate for
$q_{0}(a)$, as one can see from Example \ref{exa:cos}.
%since $u_{q}>0$ in $\overline{\Omega}$.
%\item \marginpar{\it\footnotesize added, but is it better to state Theorem \ref{c1} with in $W^{2,r}(\Omega)$ ?}Let $r>N$. Then, Theorem \ref{c1} is valid with the convergence in $W^{2,r}(\Omega)$ in place of that in $C^1(\overline{\Omega})$, as in Corollary \ref{cormt2}.
%\end{enumerate}
%\end{rem}
As a direct consequence of Theorem \ref{c1}, we complement 
%\marginpar{modified}
Theorem \ref{t10} (ii-b) showing that (A.0) is also
sufficient for the existence of a positive solution of $(P_{\mathcal{N}})$,
for some $q \in(0,1)$:

\begin{cor}
\label{cor} $(P_{\mathcal{N}})$ has a positive solution (or a solution $u\gg
0$) for \textit{some} $q \in(0,1)$ if and only if (A.0) holds.
\end{cor}

\begin{rem}
\label{rem:geN} Differently from the Dirichlet case,
%It should be mentioned that,
under (A.0) and (A.1) one may deduce the existence of a \textit{dead core}
limit function for nontrivial solutions of $(P_{\mathcal{N}})$ as
$q\rightarrow0^{+}$. Indeed, thanks to an \textit{a priori} bound
%\marginpar{modified}
\cite[Proposition 2.1]{KRQUnodea}, we may assume that a
nontrivial solution $u_{n}$ of $(P_{\mathcal{N}})$ with $q=q_{n}%
\rightarrow0^{+}$ converges to $u_{0}\geq0$ in $C^{1}(\overline{\Omega})$. We
claim that $u_{0}$ \textit{vanishes somewhere} in $\Omega$. Indeed, if
$u_{0}>0$ in $\Omega$ then Lebesgue's dominated convergence theorem shows that
$\int_{\Omega}\nabla u_{0}\nabla v=\int_{\Omega}a(x)v$ for all $v\in
C^{1}(\overline{\Omega})$, i.e. $u_{0}$ is a nontrivial solution of
$(P_{\mathcal{N}})$ with $q=0$, implying $\int_{\Omega}a=0$, a contradiction.
This situation does not occur in $(P_{\mathcal{D}})$ under (A.3$^{\prime}$)
(see Theorem \ref{cormt2}).
%\newline

\end{rem}

The final result of this section is a characterization of the set
$\mathcal{I}_{\mathcal{N}}$, which is proved by combining the \textbf{IFT} and
the sub-supersolutions method \cite[Theorem 1.4 (i)]{KRQUnodea}. Also, using
the \textbf{IFT} approach developed by Brown and Hess \cite[Theorem 1]{BH90},
we have a stability result analogous to the one in Theorem \ref{mt2c}.
%we can prove the asymptotical stability for \textit{every} solution $u\gg0$ of $(P_{\mathcal{N}})$ with $q\in \mathcal{I}_{\mathcal{N}}$ ().

\begin{theorem}
\label{th}Assume (A.0). Then $\mathcal{I}_{\mathcal{N}}=(\hat{q}_{\mathcal{N}%
},1)$ for some $\hat{q}_{\mathcal{N}}\in\lbrack0,1)$. Moreover, for
$q\in\mathcal{I}_{\mathcal{N}}$, the unique solution
%\marginpar{removed, since it is mentioned in section 7}
in $\mathcal{P}_{\mathcal{N}}^{\circ}$ is asymptotically stable.
\end{theorem}

\begin{figure}[tbh]
\begin{center}
\includegraphics[scale=0.2]{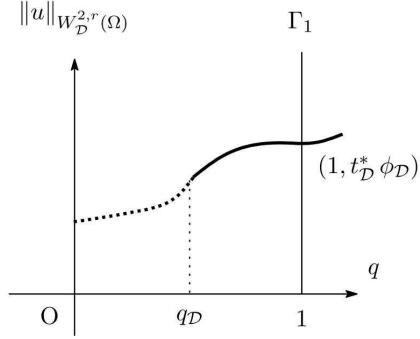}
\end{center}
\caption{The bifurcation curve of the unique nontrivial solution in the case
$\mu_{\mathcal{D}}(a)=1$, assuming that $\Omega_{+}$ is connected, satisfies
(A.2), and includes a tubular neighborhood of $\partial\Omega$. Here the full
curve represents $u_{\mathcal{D}}(q)\gg0$, whereas the dotted curve represents
solutions vanishing somewhere in $\Omega$.}%
\label{fig04}%
\end{figure}

In addition to the local result given by Theorem \ref{c1}, we can give a
\textit{global} description (i.e. for all $q\in(0,1)$) of the nontrivial
solutions set of $(P_{\mathcal{N}})$ when $\Omega_{+}$ is connected
%\marginpar{modified}
and satisfies (A.2):

\begin{rem}
\label{rem:nss} If $\Omega_{+}$ is connected and satisfies (A.2), then
Corollary \ref{old6.2} yields that $u_{\mathcal{D}}(q)\gg0$ for $q\in
(q_{\mathcal{D}},1)$, and the unique nontrivial solution of $(P_{\mathcal{D}%
})$ does not belong to $\mathcal{P}_{\mathcal{D}}^{\circ}$ for $q\in
(0,q_{\mathcal{D}}]$. Moreover, if additionally $\Omega_{+}$ includes a
tubular
%\marginpar{fixed typo}
neighborhood of $\partial\Omega$, then this solution vanishes somewhere in
$\Omega$. Note that the asymptotic behavior of $u_{\mathcal{D}}(q)$ as
$q\rightarrow1^{-}$, i.e. assertions (i) and (ii) of Theorem \ref{cormt2},
remain valid, assuming additionally (A.4), see Figure \ref{fig04}. A similar
result holds for $(P_{\mathcal{N}})$ if we assume, in addition, (A.0). In this
case, the asymptotic behavior of the solution $u_{q}\gg0$ as $q\rightarrow
1^{-}$, i.e. assertions (i)-(iii) of Theorem \ref{c1}, also remain valid
without assuming (A.4). \newline
\end{rem}

%--------
%--------

\section{Some further results}

In this section we present some results %\marginpar{modified}
(without proofs) on the two following issues:

- \textit{Explicit }sufficient conditions for the existence of positive
solutions for $(P_{\mathcal{D}})$ and $(P_{\mathcal{N}})$.

- Sufficient conditions for the existence of dead core solutions for
$(P_{\mathcal{N}})$.\smallskip

Given $0<R_{0}<R$, we write $B_{R_{0}}:=\left\{  x\in\mathbb{R}^{N}:\left\vert
x\right\vert <R_{0}\right\}  $. When $\Omega=B_{R}$ and $a$ is radial, we
shall exhibit some explicit conditions on $q$ and $a$ so that $(P_{\mathcal{D}%
})$ and $(P_{\mathcal{N}})$ admit a positive solution. In Theorem \ref{cc}
below we consider the case that $\text{supp }a^{+}$ is contained in $B_{R_{0}%
}$ and we give a condition that guarantees the existence of a positive
solution $u$ (not necessarily $\gg0$),
%\marginpar{\textit{removed, and added}}
while in Theorem \ref{rad2} we consider the case that $\text{supp }a^{-}$ is
contained in $B_{R_{0}}$ and we provide a solution $u\gg0$. These theorems are
based on a sub-supersolutions approach and are inspired in the proofs of
\cite[Section 3]{GK1} (see also the proof of Theorem \ref{t11} (iii)). If $f$
is a radial function, we shall write $f\left(  x\right)  :=f\left(  \left\vert
x\right\vert \right)  :=f\left(  r\right)  $, and we also set $A_{R_{0}%
,R}:=\left\{  x\in\mathbb{R}^{N}:R_{0}<\left\vert x\right\vert <R\right\}  $.

\begin{theorem}
\label{cc} Let $a\in C(\overline{B_{R}})$
%L^{\infty}(\Omega)$
be a radial function such that

\begin{itemize}
\item $a\geq0$ in $B_{R_{0}}$ and $a\leq0$ in $A_{R_{0},R}$;
%\marginpar{modified}

\item $r\rightarrow a(r)$ is differentiable and nonincreasing in $(R_{0},R)$,
and
\begin{equation}
\frac{1-q}{1+q}\int_{A_{R_{0},R}}a^{-}\leq\int_{B_{R_{0}}}a^{+}.
\label{inferno}%
\end{equation}
\noindent Then, $(P_{\mathcal{D}})$ has a positive solution. If, in addition,
(A.0) holds, then $(P_{\mathcal{N}})$ has a positive solution.
\end{itemize}
\end{theorem}

\noindent Note that
%\marginpar{added}
(\ref{inferno}) holds for all $q\in\left(  0,1\right)  $ if $\int_{B_{R}}%
a\geq0$, and this condition can also be formulated as
\begin{equation}
\frac{-\int_{\Omega}a}{\int_{\Omega}\left\vert a\right\vert }\leq q<1.
\label{cq}%
\end{equation}
In particular, we see that \eqref{inferno} is satisfied if $q$ is close enough
to $1$. Note that if we replace $a$ by
\[
a_{\delta}=a^{+}-\delta a^{-},\ \mbox{ with }\ \delta>\delta_{0}:=\frac
{\int_{\Omega}a^{+}}{\int_{\Omega}a^{-}},
\]
%keep $a^+$ fixed and make $a^-$ large
then the left-hand side in \eqref{cq} approaches $1$ as $\delta\rightarrow
\infty$, so that this condition becomes very restrictive for $a_{\delta}$ as
$\delta\rightarrow\infty$. On the other side, $\int_{\Omega}a_{\delta
}\rightarrow0^{-}$ as $\delta\rightarrow\delta_{0}^{+}$, so that \eqref{cq}
becomes much less constraining for $a_{\delta}$ as $\delta\rightarrow
\delta_{0}^{+}$.

We denote by $\omega_{N-1}$ the surface area of the unit sphere $\partial
B_{1}$ in $\mathbb{R}^{N}$.

\begin{theorem}
\label{rad2}

Let $a\in C(\overline{B_{R}})$
%L^{\infty}\left(  \Omega\right)  $
be a radial function satisfying (A.0). Assume that $a\geq0$ in $A_{R_{0},R}$
and
\begin{equation}
\frac{1-q}{2q+N\left(  1-q\right)  }\omega_{N-1}R_{0}^{N}\left\Vert
a^{-}\right\Vert _{C(\overline{B_{R_{0}}})}<\int_{A_{R_{0},R}}a^{+}.
\label{sipi}%
\end{equation}
Then $(P_{\mathcal{N}})$ has a solution $u\gg0$.
\end{theorem}

Unlike in Theorem \ref{cc}, we observe that no differentiability nor
monotonicity condition is imposed on $a^{-}$ in Theorem \ref{rad2}. Note again
that (\ref{sipi}) is also clearly satisfied if $q$ is close enough to $1$.
\smallskip

Finally, we consider the existence of nontrivial dead core solutions of
$(P_{\mathcal{N}})$. From \cite{BPT1, BPT2} we recall that the set
$\{x\in\Omega:u(x)=0\}$ is called the \textit{dead core} of a nontrivial
solution $u$ of $(P_{\mathcal{N}})$ if it contains an interior point.
%\marginpar{removed, since we have not proven (neither in the comment after Cor. 4.2, nor in Rem. 5.9) that $q_{\mathcal{N}}>0$ (i.e., we may have that $q_{\mathcal{N}}=0$ and so all solutions are in the cone)}
Recall that in Theorem \ref{t12} we have already given sufficient conditions
for the existence of a nontrivial solution of $(P_{\mathcal{N}})$ vanishing
somewhere in $\Omega$. We proceed now with the construction of dead cores for
%\marginpar{\textit{modified, to avoid repetition}}
solutions of $(P_{\mathcal{N}})$. To this end, let us first introduce the
following assumption:
\begin{equation}
0\leq b_{1},b_{2}\in C(\overline{\Omega})\text{\quad and\quad}\text{supp
}b_{1}\cap\{x\in\Omega:b_{2}(x)>0\}=\emptyset. \label{b1b2}%
\end{equation}
Given
%\marginpar{\textit{changed twice }$G_{\rho}$\textit{ to }$G^{\rho}$\textit{, to avoid confusion with} $\Omega_{\rho}$}
a nonempty open subset $G\subseteq\Omega$ and $\rho>0$, we set
\begin{equation}
G^{\rho}:=\left\{  x\in G:\text{dist}(x,\partial G)>\rho\right\}  .
\label{sig}%
\end{equation}
The following result is based on a comparison argument from \cite{FP84}:
%\marginpar{modified}

\begin{theorem}
\label{dc}

Let $a_{\delta}:=b_{1}-\delta b_{2}$, with $b_{1},b_{2}\not \equiv 0$
satisfying \eqref{b1b2}, and $\delta>0$. If we set $G:=\{x\in\Omega
:b_{2}(x)>0\}$ then, given $0<\overline{q}<1$ and $\rho>0$, there exists
$\delta_{0}=\delta_{0}(\rho,\overline{q})>0$ such that any nontrivial solution
of $(P_{\mathcal{N}})$ with $a=a_{\delta}$ and $q\in(0,\overline{q}]$ vanishes
in $G^{\rho}$ if $\delta\geq\delta_{0}$.
%\newline

\end{theorem}

Theorem \ref{dc} holds also for the Dirichlet problem $(P_{\mathcal{D}})$. In
particular, it complements Theorem \ref{t3a} 
%\marginpar{\it\footnotesize changed}
%\cite[Theorem 1.1]{KRQU16}
as follows: given 
$q\in(0,1)$ there exist $0<\delta_{1}<\delta_{0}$ such that every nontrivial
solution $u$ of $(P_{\mathcal{D}})$ with $a=a_{\delta}$ satisfies $u\gg0$
%$u>0$ in $\Omega$ and $\frac{\partial u}{\partial\nu}<0$ on $\partial\Omega$
for $\delta<\delta_{1}$, whereas $u$ has a nonempty dead core for
%\marginpar{removed Cor. 6.4}
$\delta>\delta_{0}$. \newline

\section{Final remarks}

\label{sec:oq}Several conditions in this paper are assumed for the sake of
presentation or technical reasons. As a matter of fact, the results in
Sections \ref{sec:P} and \ref{sec:BDN} remain true more generally for $a\in
L^{r}(\Omega)$ with $r>N$. In this situation, we assume, instead of (A.1),
that
%\marginpar{\it\footnotesize label (A.1+) removed since it is not used}%
\[
\left\{
\begin{array}
[c]{l}%
\mbox{$\Omega_{+}$ is the largest open subset of $\Omega$ where $a>0$ a.e.,}\smallskip
\\
\mbox{satisfies}\ |(\mathrm{supp}\,a^{+})\setminus\Omega_{+}%
|=0\ \mbox{and has a finite number}\smallskip\\
\mbox{of connected components,}
\end{array}
\right.
%\leqno{(\mbox{\bf A.1$^+$})}
\]
where supp is the support in the measurable sense.

It is also important to highlight that the uniqueness results in Theorem
\ref{t10} hold without assuming (A.1) and (A.2).
%\marginpar{Th. 5.1 too?}
Indeed, one may prove that the ground state solution $U_{q}$ is the only
solution of $(P_{\mathcal{D}})$ being positive in $\Omega_{+}$, and a similar
result applies to $(P_{\mathcal{N}})$ under (A.0), see \cite{KRQUpams}. A
similar situation occurs in Theorem \ref{cormt2}:
%\marginpar{modified}
without (A.1) and (A.2) the solution $u_{\mathcal{D}}(q)$ still exists for
every $q\in(0,1)$, and satisfies assertions (i)-(iv) in Theorem \ref{cormt2}
(cf.\ Remark \ref{remgs} (i)).

Also, let us mention that some of the results in this paper can be extended to
the Robin problem
%\begin{subequations}
\label{0}%
\begin{equation}%
\begin{cases}
-\Delta u=a(x)u^{q} & \mbox{in $\Omega$},\\
u\geq0 & \mbox{in $\Omega$},\\
\partial_{\nu}u=\alpha u & \mbox{on $\partial \Omega$},
\end{cases}
\label{rob}%
\end{equation}
with $\alpha\in\mathbb{R}$\thinspace. Some work in this direction has already
been done in \cite{KRQU-R1, KRQU-R2}.
%\marginpar{added ref and text}
Let us note that there are striking differences between (\ref{rob}) and the
problems considered here. For instance, under (A.0)--(A.2) and some additional
assumptions, for any $q\in\mathcal{I}_{\mathcal{N}}$ fixed, there exists some
$\overline{\alpha}>0$ such that (\ref{rob}) has \textit{exactly two} strongly
positive solutions for $\alpha\in\left(  0,\overline{\alpha}\right)  $,
\textit{one} strongly positive solution for $\alpha=\overline{\alpha}$, and
\textit{no }strongly positive solutions for $\alpha>\overline{\alpha}$
\cite[Theorem 1.3]{KRQU-R2}.

It is also worth pointing out that the positivity results in Section
\ref{sec:P} can be applied to the study of positive solutions for indefinite
concave-convex equations of the form $-\Delta u=a(x)u^{q}+b(x)u^{p}$, where
$0<1<q<p$, see \cite{KRQU16, KRQUANS}. Finally, let us mention that several
results presented here can be extended to problems involving a class of fully
nonlinear homogeneous operators \cite{dSPRQ}.\newline
%Some of them are treated in this article, whereas some other ones are left to a forthcoming paper. 

We conclude now with some interesting questions that remain open in the
context of this paper:

\begin{enumerate}
\setlength{\itemsep}{0.2cm}

\item Is the set $\mathcal{I}_{\mathcal{D}}$ connected?

\item Is there some $a$ such that $\mathcal{I}_{\mathcal{N}}=\left(
0,1\right)  $\thinspace? Let us note that we can construct a sequence
$a_{n}\in L^{\infty}(\Omega)$ such that $\mathcal{I}_{\mathcal{N}}%
(a_{n})=(q_{n},1)$ with $q_{n}\searrow0$ \cite[Remark 4.5]{KRQUnodea}.
%Is there any $a$ for which $\mathcal{I}_{\mathcal{N}}(a)=(0,1)$\thinspace? Similarly, is there any $a$ for which $\mathcal{A}_{\mathcal{N}}(a)=(0,1)$\thinspace?

\item Assume $\mathcal{I}_{\mathcal{N}}(a)=(0,1)$. Can we characterize the
limiting behavior of the solution $u_{q}\gg0$ of $(P_{\mathcal{N}})$ as
$q\rightarrow0^{+}$\thinspace?

\item By Theorem \ref{t12}, we see that we may have $q_{\mathcal{D}}>0$ or
$q_{\mathcal{N}}>0$. On the other side, Theorem
%\marginpar{modified}
\ref{cormt2} (iv-c) shows a situation in which $q_{\mathcal{D}}=0$. Can we
have $q_{\mathcal{N}}=0$ (i.e., $\mathcal{A}_{\mathcal{N}}=\left(  0,1\right)
$) ?

\item Can we obtain explicit sufficient conditions for the existence of
positive solutions of $(P_{\mathcal{D}})$ and $\left(  P_{\mathcal{N}}\right)
$ (as e.g. the ones in Theorems \ref{cc} and \ref{rad2} for $\left(
P_{\mathcal{N}}\right)  $; or the ones in Theorem \ref{cc} and \cite[Theorem
3.2 (i)]{GK1} for $\left(  P_{\mathcal{D}}\right)  $) without assuming that
$\Omega$ is a ball and $a$ is radial?

\item Is it possible to extend the results in this paper to a general operator
of the form%
\[
Lu=-\operatorname{div}(A(x)\nabla u)+\langle b(x),\nabla u\rangle+c\left(
x\right)  u,
\]
under suitable assumptions on the coefficients? Let us note that if
$b\not \equiv 0$ variational techniques do not apply. Furthermore, the size of
the coefficient $c$ plays an important role: in the one-dimensional Dirichlet
case no positive solutions exist if $c>0$ is large enough, cf. \cite[Theorem
3.11]{ejde}. Let %\marginpar{\textit{added}}
us add that the Neumann case with
$A\equiv1$, $b\equiv0$, and $c$ constant has been treated in \cite{alama,
KRQU-R2}. 

\item We
%\marginpar{moved as item}
believe that many of the resuts and techniques reviewed here also apply to the
corresponding $p$-Laplacian equation
\[
-\Delta_{p}u=a(x)u^{q},
\]
with $p>1$ and $0<q<p-1$. Some progress in this direction has been achieved in
\cite{KRQUpams}.
\end{enumerate}

%The rest of the article is organized as follows. In Section \ref{secneg} we mainly analyze the case $\alpha\leq0$ and prove Theorems \ref{tpos} and \ref{thm:curve}. Section \ref{secpos} is mostly devoted to $(R_{\alpha})$ with $\alpha>0$, where we investigate qualitative properties of the solutions set and prove an exact multiplicity result employing the change of variables \eqref{changeofv}. Lastly, Section \ref{sec:topoba} provides a topological bifurcation approach of $(R_{\alpha})$ and the proof of Theorem \ref{thm:main}.

%-------------------------------------------------

%\bibitem {alama}{\small S. Alama, \textit{Semilinear elliptic equations with
%sublinear indefinite nonlinearities}, Adv. Differential Equations \textbf{4}
%(1999), 813--842. }

%\end{subequations}

\end{document}